\RequirePackage{fix-cm}
\documentclass[journal]{IEEEtran}
\usepackage{cite}
\usepackage{amsmath,amssymb,amsfonts}
\usepackage{algorithmic}
\usepackage{graphicx}
\usepackage{textcomp}
\usepackage{diagbox} 
\usepackage[utf8]{inputenc}
\usepackage{subfigure}
\usepackage{placeins}
\usepackage{amssymb}
\usepackage[numbers,sort&compress]{natbib}

\usepackage{caption}
\usepackage{tabularx}
\usepackage{subfigure}
\usepackage{multirow}
\usepackage{booktabs}
\usepackage{graphicx}
\usepackage{lipsum}
\usepackage{multicol}
\usepackage{listings}
\usepackage{xcolor} 
\newlength\savedwidth

\newlength\savewidth

\usepackage{soul}
\usepackage{soul}
\soulregister\cite7 
\soulregister\citep7 
\soulregister\citet7 
\soulregister\ref7 
\soulregister\pageref7 

\begin{document}
	
	\title{Mathematical representation of the WECC composite load model}

	
\author{Zixiao~Ma,~\IEEEmembership{Student Member,~IEEE,}~
	Zhaoyu~Wang,~\IEEEmembership{Member,~IEEE,}~
	Yishen~Wang,~\IEEEmembership{Member,~IEEE,}~
	Ruisheng~Diao,~\IEEEmembership{Senior~Member,~IEEE,}~and~
	Di~Shi,~\IEEEmembership{Senior~Member,~IEEE}
	\thanks{Zixiao Ma and Zhaoyu Wang are with the Department of Electrical and Computer Engineering, Iowa State University, Ames, IA, 50011 USA e-mail: zma@iastate.edu; wzy@iastate.edu).}
	\thanks{Yishen Wang, Ruisheng Diao and Di Shi are with GEIRI North America, San Jose, CA 95134, USA (email: yishen.wang@geirina.net, ruisheng.diao@geirina.net, di.shi@geirina.net).}
	
}
	

	\maketitle
	
	\begin{abstract}
		The Western Electricity Coordinating Council (WECC) composite load model is a newly developed load model that has drawn great interest from the industry. To analyze its dynamic characteristics with both mathematical and engineering rigor, a detailed mathematical model is needed. Although WECC composite load model is available in commercial software as a module and its detailed block diagrams can be found in several public reports, there is no complete mathematical representation of the full model in literature. This paper addresses a challenging problem of deriving detailed mathematical representation of WECC composite load model from its block diagrams. In particular, for the first time, we have derived the mathematical representation of the new DER\_A model. The developed mathematical model is verified using both Matlab and PSS/E to show its effectiveness in representing WECC composite load model. The derived mathematical representation serves as an important foundation for parameter identification, order reduction and other dynamic analysis. 
	\end{abstract}
	\begin{IEEEkeywords}
		Composite load model, dynamic load modeling, mathematical model, three-phase motor, DER\_A.
	\end{IEEEkeywords}
	\section{Introduction}
	\label{sec:introduction}
	Load modeling is essential to power system stability analysis, optimization, and controller design as shown in many research \cite{C.W.Taylor1994}. Although the importance of load modeling is recognized by power system researchers and engineers \cite{P.Kundur1994}, obtaining an accurate load model remains challenging. The difficulty is caused by the large number of diverse load components, time-varying compositions, and the lack of detailed load information and measurements. To this end, developing high-fidelity load models that approximate the real load characteristic while overcoming the above challenges is imperative.
	
	Load modeling consists of developing model structures and identifying associated parameters. For a given load model structure, its parameter identification can be implemented with component or measurement-based approaches. The component-based approach is based on the knowledge of detailed physical models of different load components and their compositions. \cite{Wong2012,Kosterev2008}. However, such information is usually difficult to obtain, which motivates the research of measurement-based load modeling \cite{LeeJune2012,Ma2008,Visconti2014,Han2009,Choi2009,Hu2016}. With the wider deployment of digital fault recorders, the measurement-based load modeling approaches become increasingly popular \cite{Choi2009,Ma2008,Son2014,Kim2016,Overbye2012}. Measurement-based modeling uses the measured data to identify model parameters. The main advantage of this approach is that it collects the data directly from the power system and can be used for online modeling.
	
\begin{figure}[ht]
		\centering
		\includegraphics[width=8cm]{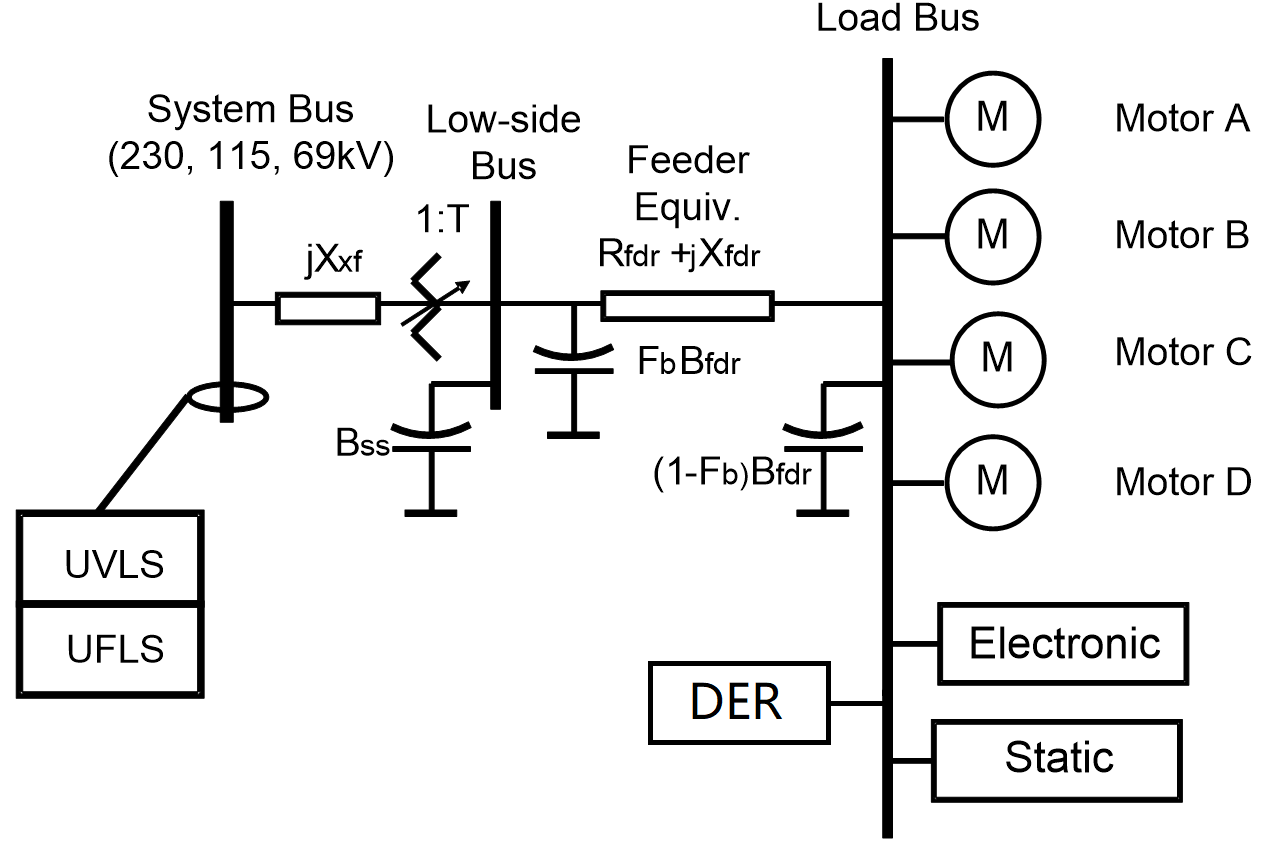}
		\caption{A schematic of the WECC CMPLDWG composite load model \cite{WECC2015}.}
		\label{CMPLDWG}
	\end{figure}

	For the load model structures, there exists static and dynamical load models. For example, static load models includes static constant impedance-current-power (ZIP) model and exponential model \cite{Kosterev2008}. However, they cannot capture the dynamic behaviors of loads. Dynamic load models represent the real/reactive powers as functions of both voltage and time, such as Induction Motor (IM) model and Exponential Recovery Load Model (ERL) \cite{Hiskens2001,Arif2018,Huang2019}. To consider both dynamic and static load characteristics, composite load models are proposed, such as ZIP+IM load model, Complex Load Model (CLOD), Low-Voltage (LV) Load Models and WECC composite load model (WECC CLM). An aggregated five-machine dynamic equivalent electro-mechanical model of the WECC power system using synchrophasor measurements was developed to bridge the gap aroused by the increasing penetration of renewable energy resources. These renewable resources will significantly change dynamic properties, inter-area oscillation characteristics and stability margins of WECC power systems in the near future \cite{Chavan2016}. However, this model is built from the entire power system's point of view. After the 1996 blackout of the Western Systems Coordinating Council (WSCC) \cite{Kosterev1999}, the ZIP+IM model was designed to capture the dynamic effects of highly stressed conditions in summer peaks. However, this interim model was ineffective in capturing delayed voltage recovery events from transmission faults \cite{Kosterev2008,Williams1992,Shaffer1997}. By adding the electrical distance between the transmission system and the electrical end-uses, as well as adding special components such as electronic load components and single-phase motors, a preliminary WECC CLM was proposed and implemented in major industry-level commercial  simulation software packages \cite{Arif2018}. With continuous updates and the incorporation of distributed energy resources (DERs), the newest developed WECC composite load model called CMPLDWG is proposed as shown in Figure \ref{CMPLDWG}. The model includes an electrical representation of a distribution system with a substation transformer, shunt reactance, and a feeder equivalent. At distribution system side, it includes a static load model, one power electronics model, three three-phase motor models, one AC single phase motor and a distributed energy resource. CMPLDWG uses PVD1 model to represent the DERs. However, PVD1 constitutes a total of 5 modules, 121 parameters and 16 states , which is as complex as the WECC CLM. Therefore, EPRI has developed a simpler yet more comprehensive model to replace PVD1, which is named as DER\_A model \cite{EPRI2019}.
	
	Although the WECC composite load model has been widely implemented in commercial power system software, a comprehensive mathematical representation cannot be found in existing literature. Moreover, researchers cannot access the source codes of commercial software packages, making it hard to obtain insights of the models implemented in the software. The detailed block diagrams of the model can be found in publicly available reports, for example \cite{Wang2018.,EPRI2019}. However, deriving mathematical representation from these diagrams are challenging. In \cite{Huang2017a}, a mathematical representation of three-phased motors has been provided, nevertheless, the DER\_A model is missing. However, the mathematical model is essential for parameter identification, stability assessment, and dynamic order reduction. To this end, this paper derives a detailed and comprehensive mathematical representation of the WECC composite load model with DER\_A. Various simulations are conducted in both matlab and PSS/E to verify the effectiveness of the derived mathematical model. 
	
	The rest of the paper is organized as follows. Section \ref{C2} presents the detailed derivation of mathematical model of WECC composite load model. Section \ref{C3} shows the simulation results and analysis. Section \ref{C4} concludes the paper.
	\section{Mathematical modeling of individual components}\label{C2}
	In this section, we will derive mathematical representations for individual components in WECC composite load model, namely, three-phase motors, DER\_A, single-phase motor, electronic and static loads. 
	\subsection{Three-phase motor model}
	There are multiple types of three-phase induction motors that can describe the end-use loads \cite{NARC1}. In WECC CMPLDWG three different three-phase motors, A, B and C are used to represent different types of dynamic components. Motor A represents three-phase induction motors with low inertia driving constant torque loads, e.g. air conditioning compressor motors and positive displacement pumps. Motor B represents three-phase induction motors with high inertia driving variable torque loads such as commercial ventilation fans and air handling systems. Motor C represents three-phase induction motors with low inertia driving variable torque loads such as the common centrifugal pumps. 
	
\begin{figure}[ht]
		\centering
		\includegraphics[width=8.5cm]{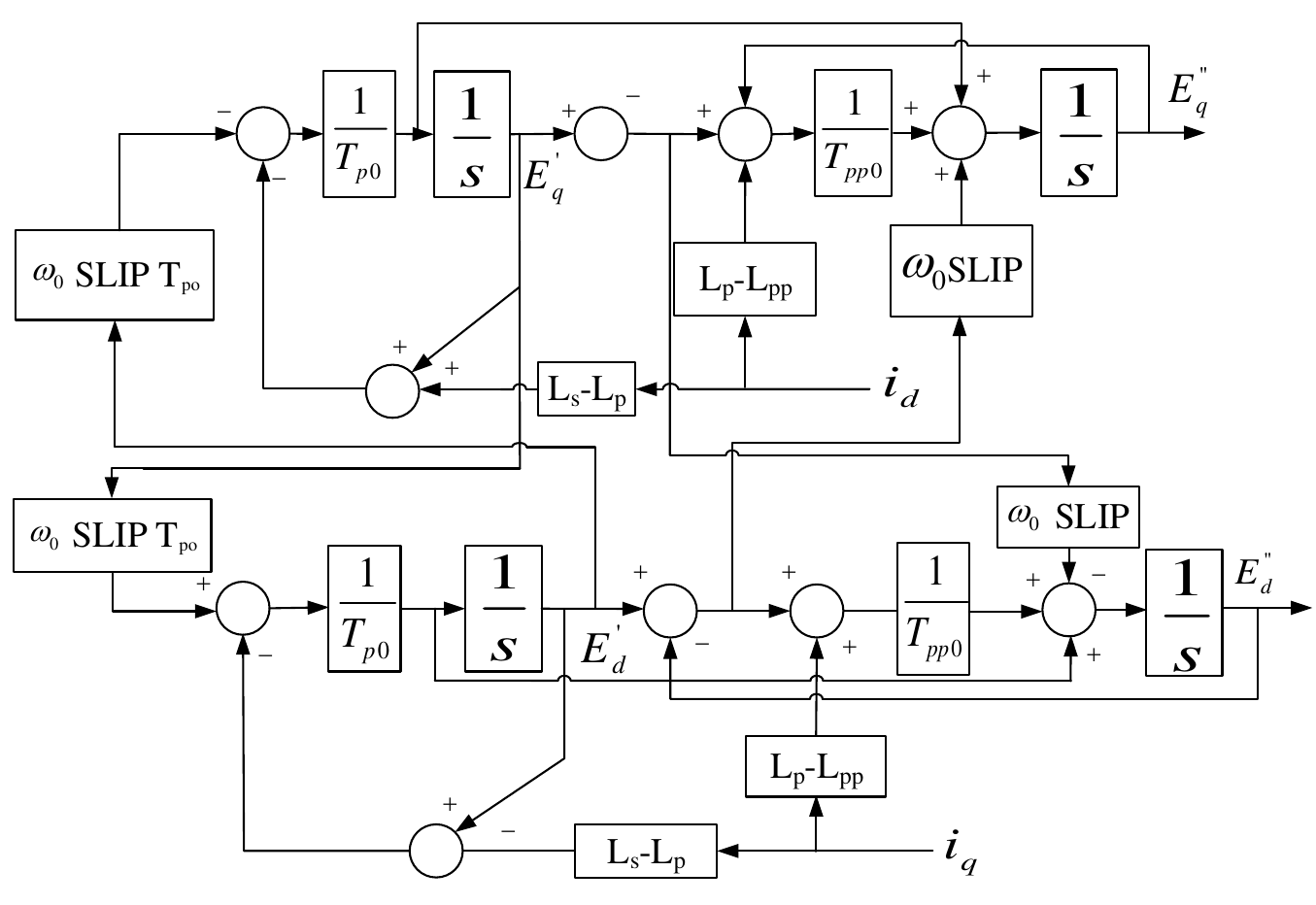}
		\caption{The diagram of three-phase motor.}
		\label{threephase}
	\end{figure}

	These three-phase motors share the same model structure, however, their model parameters are different. Therefore, a fifth-order induction motor model is adopted to represent three-phase motors in the WECC composite load model. The block diagram is shown in Figure \ref{threephase}. From the diagram we can obtain a fourth-order electrical model with respect to ${E}_q'$, ${E}_d'$, ${E}_q{''}$ and ${E}_d{''}$. Combining with the mechanical model, we have the complete fifth-order model as follows,
	\begin{align}
	&\dot {E}_q'\! = \!\frac{1}{{{T_{p0}}}}\left[ { - E_q'\! -\! {i_d}\left( {{L_s} \!- \!{L_p}} \right)\! - \!E_d'\! \cdot \!{\omega _0} \!\cdot\! SLIP \!\cdot\! {T_{P0}}} \right],\\\label{Eq'}
	&\dot {E}_d' \!=\! \frac{1}{{{T_{p0}}}}\left[ {-E_d'\! + \!{i_q}\left( {{L_s}\! -\! {L_p}} \right) \!+ \!E_q' \!\cdot\! {\omega _0} \!\cdot\! SLIP\! \cdot\! {T_{P0}}} \right],\\\label{Ed'}
	&\dot E_q{''} \!\!={\frac{{T_{p0}\!-\!T_{pp0}}}{T_{p0}T_{pp0}}}E_q'\!-\!{\frac{T_{pp0}\left({ L_s\!\!-\!\!L_p}\right) \!+\!T_{p0}\left( {L_p\!\!-\!\!L_{pp}}\right) }{{T_{p0}T_{pp0}}}}i_d\!\!\\\nonumber
	&	-{\frac{1}{T_{pp0}}}E_q''-\omega_0\cdot SLIP \cdot E_d'',\\\label{Eq''}
	&\dot E_d{''} \!\!={\frac{{T_{p0}\!-\!T_{pp0}}}{T_{p0}T_{pp0}}}E_d'\!+\!{\frac{T_{pp0}\left({ L_s\!\!-\!\!L_p}\right) \!+\!T_{p0}\left( {L_p\!\!-\!\!L_{pp}}\right) }{{T_{p0}T_{pp0}}}}i_q\!\!\\\nonumber
	&-{\frac{1}{T_{pp0}}}E_d''+\omega_0\cdot SLIP \cdot E_q'',\\\label{Ed''}
	&S\dot LIP =  - \frac{{p \cdot E_d{''} \cdot {i_d} + q \cdot E_q{''} \cdot {i_q} -TL}}{{2H}}.
	\end{align}
	The algebraic equations are:
	\begin{align}
	&TL = {T_{m0}}\left( {A{w^2} + Bw + {C_0} + D{w^{Etrq}}} \right),\\
	&{T_{m0}} = pE_{d0}{''}{i_{d0}} + qE_{q0}{''}{i_{q0}},\\
	&w = 1 - SLIP,\\
	&i_{d}={\frac{r_s}{r_s^2+L_{pp}^2}}({V_d+E_d{''}} )+ {\frac{L_{pp}}{r_s^2+L_{pp}^2}}({V_q+E_q{''}} ), \\
	&i_{q}={\frac{r_s}{r_s^2+L_{pp}^2}}({V_q+E_q{''}} )- {\frac{L_{pp}}{r_s^2+L_{pp}^2}}({V_d+E_d{''}}), \\
	&P=V_{d}i_d+V_{q}i_q, \\
	&Q=V_{q}i_d-V_{d}i_q, 
	\end{align}
	where the five state variables are ${E}_q'$, ${E}_d'$, ${E}_q{''}$, ${E}_d{''}$ and $SLIP$; $L_s$, $L_p$ and $L_{pp}$ are synchronous reactance, transient and subtransient
	reactance, respectively; $T_{p0}$ and $T_{pp0}$ are transient and subtransient rotor time constants, respectively; and $\omega_0$ is the synchronous frequency.
	\subsection{Single-phase motor model}
	Motor D in Figure \ref{CMPLDWG} represents the single-phase motor model that captures behaviors of single-phase air Figure with reciprocating compressors. However, it is challenging to model the fault point on wave \cite{Huang2016} and voltage ramping effects \cite{NARC1}. Moreover, new A/C motors are mostly equipped with scroll compressors and/or power electronic drives, making their dynamic characteristics significantly different than conventional motors. Therefore, WECC uses a performance-based model to represent single-phase motors. The main purpose of deriving the mathematical model is to establish the foundation for theoretical studies such as parameter identification and order reduction. Hence, for this purpose, it is unnecessary to derive the mathematical representation of the performance model.
	
	\begin{figure}[ht!]
		\centering
		\includegraphics[width=8cm]{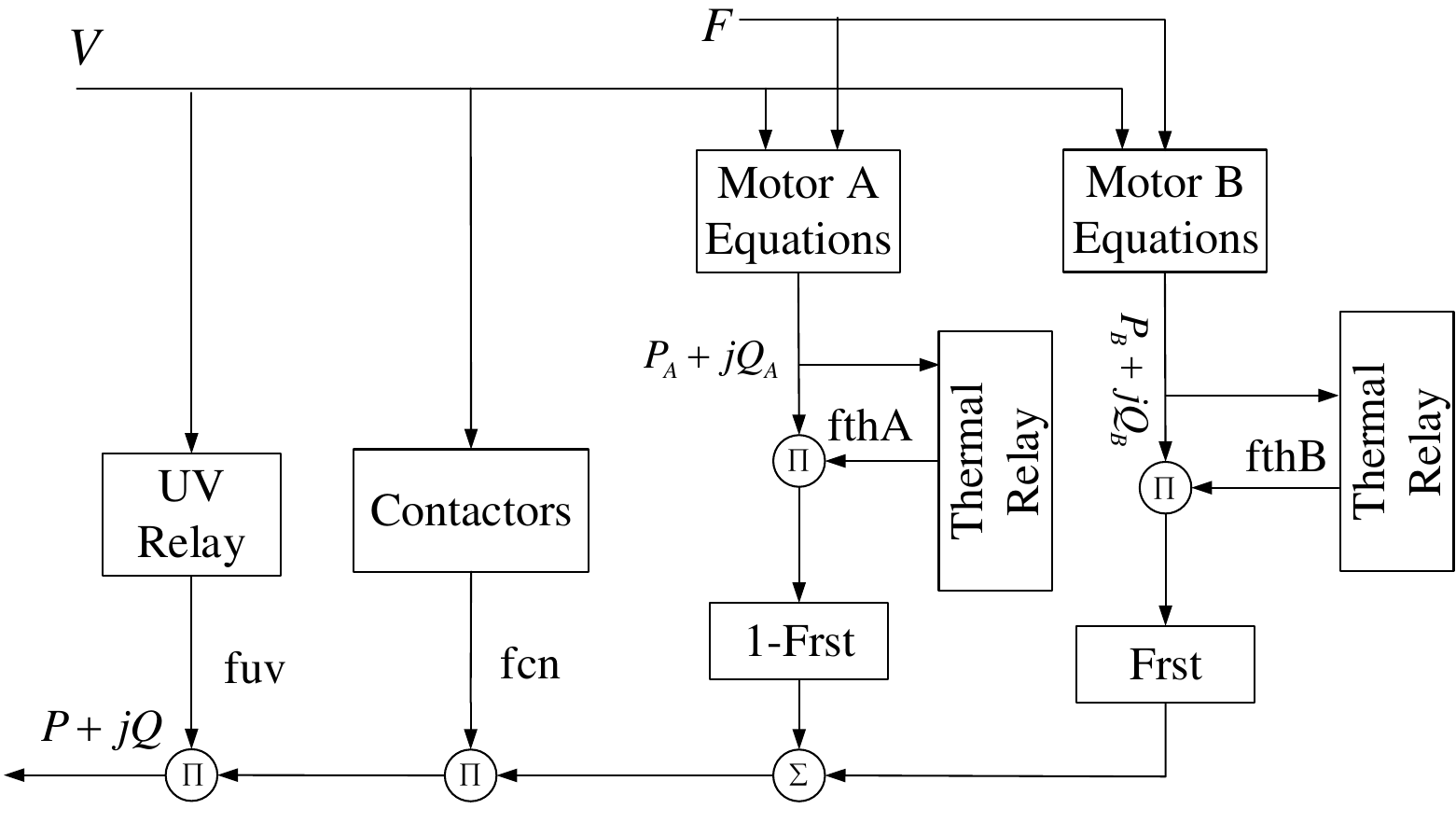}
		\caption{The diagram of single-phase motor \cite{Huang2017a, WECC2015}.}
		\label{singlephase}
	\end{figure}
	\subsection{DER\_A model}
\begin{figure*}[ht!]
		\centering
		\includegraphics[width=17.5cm]{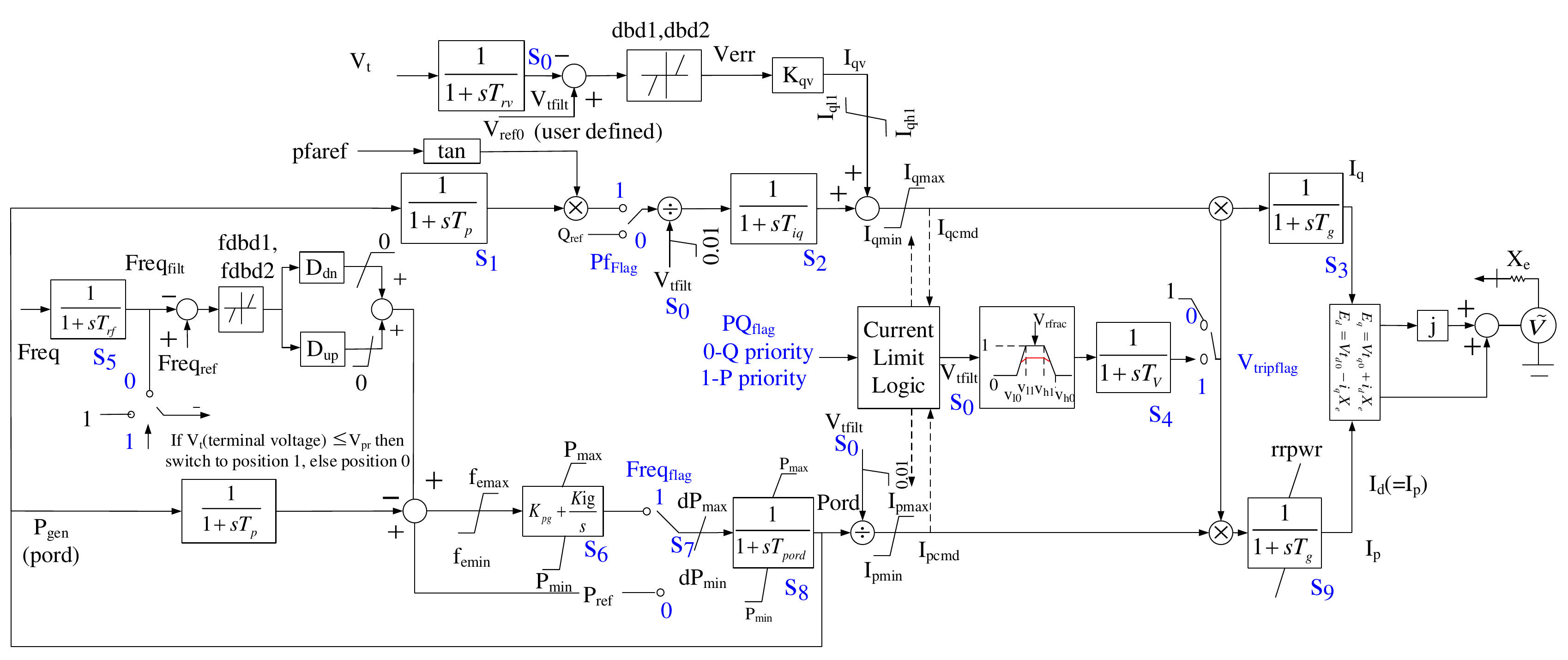}
		\caption{The diagram of DER\_A model \cite{EPRI2019}.}.
		\label{DERA}
	\end{figure*}

	The DER\_A is a newly developed model representing aggregate renewable energy resources. Compared to the previous PVD1 model that is relatively large scale and complex, the DER\_A model has fewer states and parameters. There is no mathematical representation of the DER\_A model in existing literature till now. In this section, the detailed mathematical model is derived from Figure \ref{DERA} with respect to each state variable. The parameters are defined in Table \ref{table1}.
	\subsubsection{Mathematical model of $S_0$}
	Figure \ref{1} shows the block diagram of first-order filter whose input is the bus voltage $V_t$, and the output is filtered voltage $S_0$ ($V_{t\_filt}$). From the diagram, we can obtain the following dynamic equation:
	\begin{figure}[ht!]
		\centering
		\includegraphics[width=7cm]{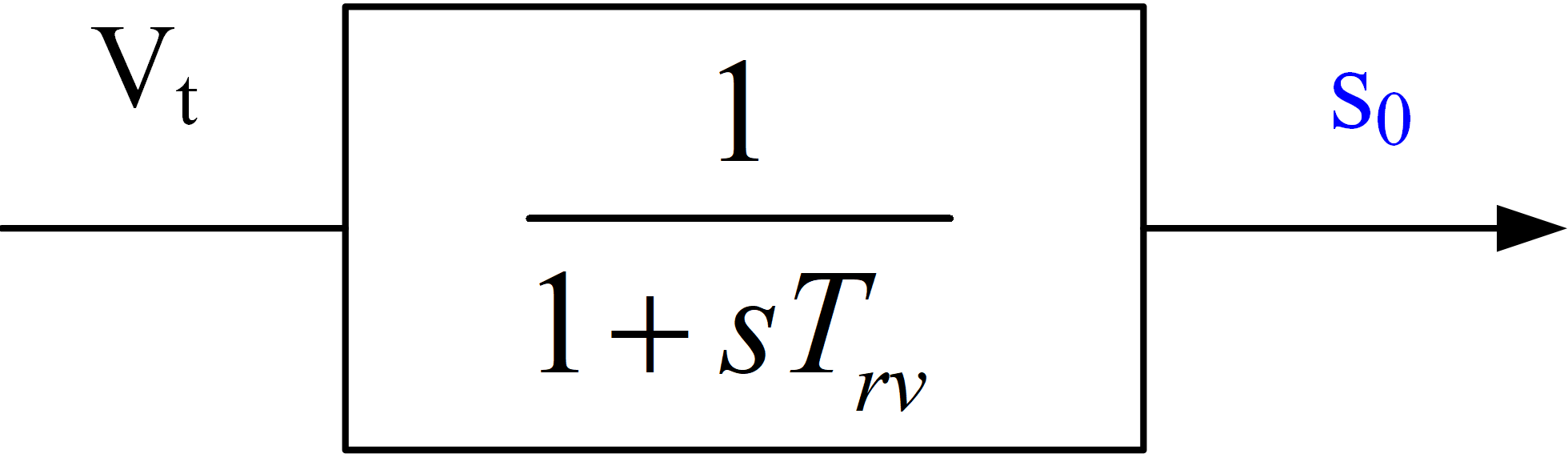}
		\caption{A local schematic of $S_0$ in the DER\_A model.}
		\label{1}
	\end{figure}
	\begin{equation}
	{\dot S_0} = \frac{1}{{{T_{rv}}}}\left( {{V_t} - {S_0}} \right)
	\end{equation}
	\subsubsection{Mathematical model of $S_1$}
	Figure \ref{2} shows the block diagram of first-order filter whose input is the electrical power being generated at the terminals of the DER\_A model $P_{gen}$, and the output is filtered power $S_1$. From the diagram, we can obtain the following dynamic equation:
	\begin{figure}[ht!]
		\centering
		\includegraphics[width=7cm]{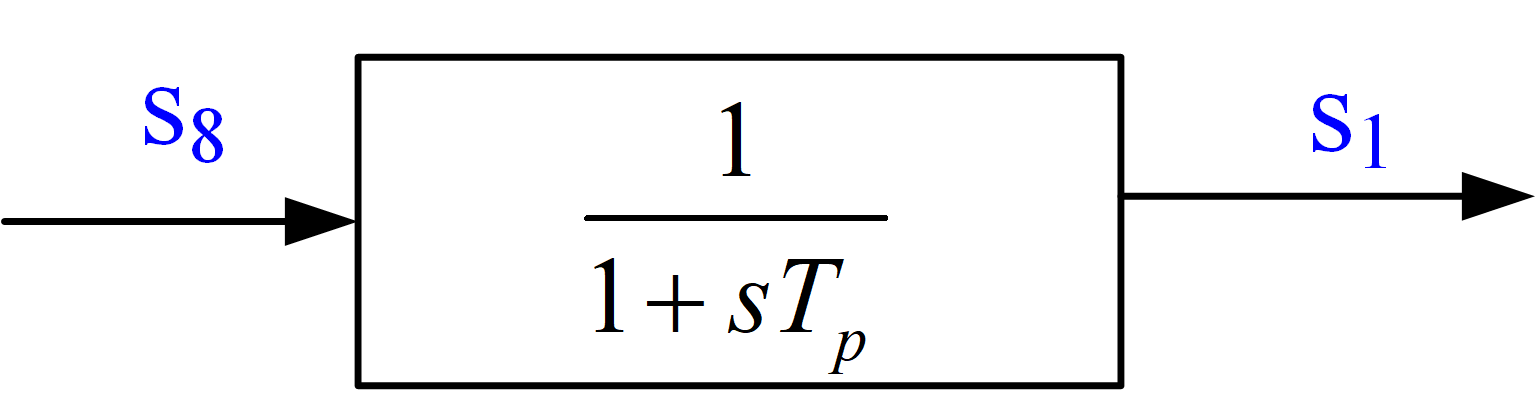}
		\caption{A local schematic of $S_1$ in the DER\_A model.}
		\label{2}
	\end{figure}
	\begin{eqnarray}
	{\dot S_1} = \frac{1}{{{T_p}}}\left( {{S_{8}} - {S_1}} \right)
	\end{eqnarray}
	\subsubsection{Mathematical model of $S_2$}
	The local block diagram of $S_2$ is shown in Figure \ref{3}. From the diagram, we can obtain the following dynamic equation:
	\begin{figure}[ht!]
		\centering
		\includegraphics[width=8cm]{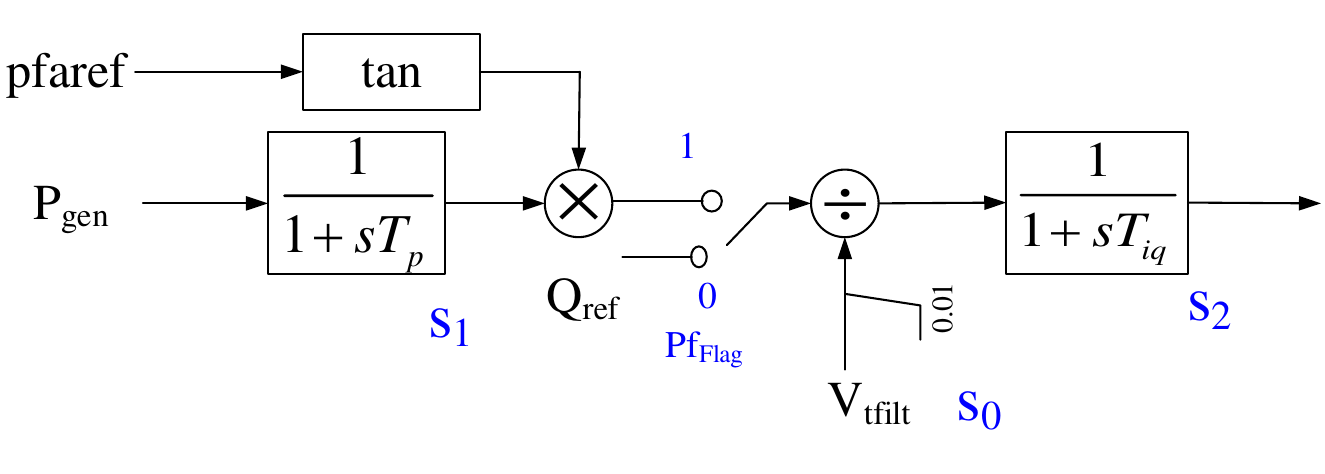}
		\caption{A local schematic of $S_2$ in the DER\_A model.}
		\label{3}
	\end{figure}
	\begin{eqnarray}
	{\dot S_2} = \left\{ \begin{aligned}
	&\!\!\!-\frac{{S_2}}{{T_{iq}}}\!  + \!\frac{{{Q_{ref}}}}{{T_{iq}}\cdot{sat_1\left( {{S_0}} \right)}} {\kern 20pt} if{\kern 1pt} {\kern 1pt} {\kern 1pt} {P_{fFlag}} = 0\\
	&\!\!\!-\frac{S_2}{{T_{iq}}} \!+ \!\frac{{\tan \left( {pfaref} \right) \!\times\! {S_1}}}{{T_{iq}}\cdot{sa{t_1}\left( {{S_0}} \right)}} {\kern 4pt} if{\kern 1pt} {\kern 1pt} {\kern 1pt} {\kern 1pt} {P_{fFlag}} = 1
	\end{aligned} \right.
	\end{eqnarray}
	where $Q_{ref}$ ($P_{ref}$ in Equation (\ref{s6})) is determined based on the initial P/Q output of the DER\_A model in software; $pfaref$  can be computed by $arctan\left(Q_{gen0}/P_{gen0}\right)$, where $Q_{gen0}$ and
	$P_{gen0}$ are the active and reactive power determined by the initial power
	flow solution. The limiter in the diagram is described by a saturation function  that can be defined as Equation (\ref{sat1}).
	\begin{equation}\label{sat1}
	sa{t_1}\left( {{x}} \right) = \left\{ \begin{aligned}
	{x}{\kern 1pt} {\kern 1pt} {\kern 1pt} {\kern 1pt} {\kern 1pt} {\kern 1pt} {\kern 1pt} {\kern 1pt} {\kern 1pt} {\kern 1pt} {\kern 1pt} {\kern 1pt} {\kern 1pt} {\kern 1pt} {\kern 1pt} {\kern 1pt} {\kern 1pt} {\kern 1pt} {\kern 1pt} {\kern 1pt} {\kern 1pt} {\kern 1pt} {\kern 1pt} {\kern 1pt} {\kern 1pt} {\kern 1pt} {\kern 1pt} {\kern 1pt} {\kern 1pt} {\kern 1pt} if{\kern 1pt} {\kern 1pt} {\kern 1pt} {x} \geqslant 0.01\\
	0.01{\kern 1pt} {\kern 1pt} {\kern 1pt} {\kern 1pt} {\kern 1pt} {\kern 1pt} {\kern 1pt} {\kern 1pt} {\kern 1pt} {\kern 1pt} {\kern 1pt} {\kern 1pt} {\kern 1pt} {\kern 1pt} {\kern 1pt} {\kern 1pt} {\kern 1pt} {\kern 1pt} {\kern 1pt} {\kern 1pt} if{\kern 1pt} {\kern 1pt} {\kern 1pt} {\kern 1pt} {x} \leqslant 0.01
	\end{aligned} \right.{\kern 1pt} {\kern 1pt} {\kern 1pt} {\kern 1pt} {\kern 1pt} {\kern 1pt} {\kern 1pt} {\kern 1pt} {\kern 1pt} {\kern 1pt}
	\end{equation}
	\subsubsection{Mathematical model of $S_3$}
	The local block diagram of q-axis current $S_3$ ($i_q$) is shown in Figure \ref{4}. From the diagram, we can obtain the following dynamic equation:
	\begin{figure}[ht!]
		\centering
		\includegraphics[width=8cm]{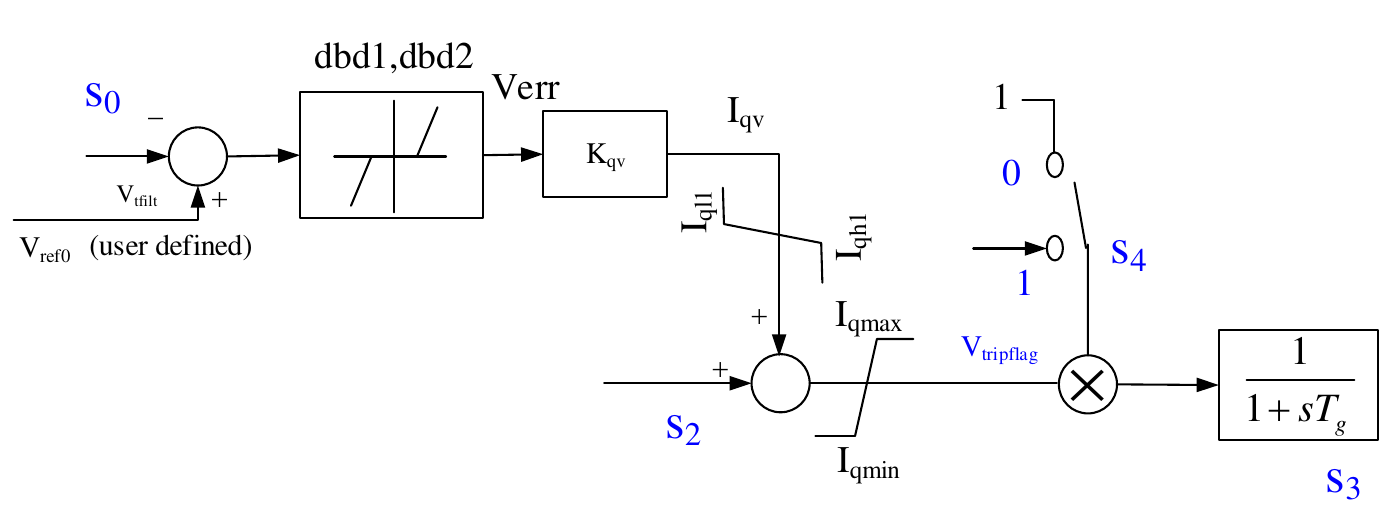}
		\caption{A local schematic of $S_3$ in the DER\_A model.}
		\label{4}
	\end{figure}
	\begin{eqnarray}
	{\dot S_3}\! \!= \!\!\left\{\! \begin{aligned}
	&\!\!- \!\!\frac{{S_3}\! -\! sat_2\left( {{S_2} \!+ \!sat_3\left( {D{B_V}\left( {{V_{ref0}} \!-\! {S_0}} \right)\! \cdot\! {K_{qv}}} \right)\! } \right)}{T_g}{\kern 1pt}  \!\!\!\!\!\!\!\!\!\!\!\!\!\!\!\!\!\!\!\!\!\!\!\!\!\\\!\!\!\!\!\!\!\!\!\!\!\!\!\!\!\!\!\!\!\!
	&{\kern 120pt}if{\kern 1pt} {\kern 1pt} {\kern 1pt} {V_{tripFlag}} = 0\!\!\!\!\!\!\!\!\!\!\!\!\!\!\!\!\!\!\!\!\\\!\!\!\!\!\!\!\! \\\!\!\!\!\!\!\!\!\!\!\!\!\!\!\!\!\!\!
	&\!\!-\!\! \frac{{S_3} \!-\! sat_2\left( {{S_2}\! +\! sat_3\left(\! {{D{B_V}\left( {{V_{ref0}} \!- \!{S_0}} \right)\! \cdot \!{K_{qv}}}} \right)\! } \right)\!\! \times\!\! {S_4}}{T_g}{\kern 1pt} {\kern 1pt} {\kern 1pt} {\kern 1pt} {\kern 1pt} {\kern 1pt} {\kern 1pt} {\kern 1pt} {\kern 1pt} {\kern 1pt} {\kern 1pt} {\kern 1pt} \!\!\!\!\!\!\!\!\!\!\!\!\!\!\!\!\!\!\!\!\!\!\!\\\!\!\!\!\!\!\!\!\!\!\!\!\!\!\!\!\!\!\!
	&{\kern 120pt}if{\kern 1pt} {\kern 1pt} {\kern 1pt} {\kern 1pt} {V_{tripFlag}} = 1\!\!\!\!\!\!\!\!\!\!\!\!\!\!\!\!\!\!\!\!\!\!\!\!\!\!
	\end{aligned} \right.{\kern 1pt} {\kern 1pt} {\kern 1pt} {\kern 1pt} {\kern 1pt} {\kern 1pt} {\kern 1pt} {\kern 1pt} {\kern 1pt} {\kern 1pt}
	\end{eqnarray}
	where the limiter function and dead bands function are defined as Eqn (\ref{sat2}) and (\ref{DBV}), respectively.
	\begin{eqnarray}\label{sat2}
	sa{t_2}\left( x \right) = \left\{ \begin{aligned}
	&{I_{q\max }}{\kern 1pt} {\kern 1pt} {\kern 1pt} {\kern 1pt} {\kern 1pt} {\kern 1pt} {\kern 1pt} {\kern 1pt} {\kern 1pt} {\kern 1pt} {\kern 1pt} {\kern 1pt} {\kern 1pt} {\kern 1pt} {\kern 1pt} {\kern 1pt} {\kern 1pt} {\kern 1pt} {\kern 1pt} {\kern 1pt} {\kern 1pt} {\kern 1pt} {\kern 1pt} {\kern 1pt} {\kern 1pt} {\kern 1pt} {\kern 1pt} {\kern 1pt} {\kern 1pt} {\kern 1pt} {\kern 1pt} {\kern 1pt} {\kern 1pt} {\kern 1pt} {\kern 1pt} {\kern 1pt} {\kern 1pt} {\kern 1pt} if{\kern 1pt} {\kern 1pt} {\kern 1pt} x \geqslant {I_{q\max }}\\ 
	&x{\kern 1pt} {\kern 1pt} {\kern 1pt} {\kern 1pt} {\kern 1pt} {\kern 1pt} {\kern 1pt} {\kern 1pt} {\kern 1pt} {\kern 1pt} {\kern 1pt} {\kern 1pt} {\kern 1pt} {\kern 1pt} {\kern 1pt} {\kern 1pt} {\kern 1pt} {\kern 1pt} {\kern 1pt} {\kern 1pt} {\kern 1pt} {\kern 1pt} {\kern 1pt} {\kern 1pt} {\kern 1pt} {\kern 1pt} {\kern 1pt} {\kern 1pt} {\kern 1pt} {\kern 1pt} {\kern 1pt} {\kern 1pt} {\kern 1pt} {\kern 1pt} {\kern 1pt} {\kern 1pt} {\kern 1pt} {\kern 1pt} if{\kern 1pt} {\kern 1pt} {\kern 1pt} {\kern 1pt} {I_{q\min }} \leqslant x \leqslant {I_{q\max }}\\ 
	&{I_{q\min }}{\kern 1pt} {\kern 1pt} {\kern 1pt} {\kern 1pt} {\kern 1pt} {\kern 1pt} {\kern 1pt} {\kern 1pt} {\kern 1pt} {\kern 1pt} {\kern 1pt} {\kern 1pt} {\kern 1pt} {\kern 1pt} {\kern 1pt} {\kern 1pt} {\kern 1pt} {\kern 1pt} {\kern 1pt} {\kern 1pt} {\kern 1pt} {\kern 1pt} {\kern 1pt} {\kern 1pt} {\kern 1pt} {\kern 1pt} {\kern 1pt} {\kern 1pt} {\kern 1pt} {\kern 1pt} {\kern 1pt} {\kern 1pt} {\kern 1pt} {\kern 1pt} {\kern 1pt} {\kern 1pt} {\kern 1pt} {\kern 1pt} {\kern 1pt} {\kern 1pt} if{\kern 1pt} {\kern 1pt} {\kern 1pt} {\kern 1pt} x \leqslant {I_{q\min }}
	\end{aligned} \right.{\kern 1pt} {\kern 1pt} {\kern 1pt} {\kern 1pt} {\kern 1pt} {\kern 1pt} {\kern 1pt} {\kern 1pt} {\kern 1pt} {\kern 1pt}
	\end{eqnarray}
	\begin{eqnarray}\label{sat3}
	sa{t_3}\left( x \right) = \left\{ \begin{aligned}
	&{I_{qh1 }}{\kern 1pt} {\kern 1pt} {\kern 1pt} {\kern 1pt} {\kern 1pt} {\kern 1pt} {\kern 1pt} {\kern 1pt} {\kern 1pt} {\kern 1pt} {\kern 1pt} {\kern 1pt} {\kern 1pt} {\kern 1pt} {\kern 1pt} {\kern 1pt} {\kern 1pt} {\kern 1pt} {\kern 1pt} {\kern 1pt} {\kern 1pt} {\kern 1pt} {\kern 1pt} {\kern 1pt} {\kern 1pt} {\kern 1pt} {\kern 1pt} {\kern 1pt} {\kern 1pt} {\kern 1pt} {\kern 1pt} {\kern 1pt} {\kern 1pt} {\kern 1pt} {\kern 1pt} {\kern 1pt} {\kern 1pt} {\kern 1pt} if{\kern 1pt} {\kern 1pt} {\kern 1pt} x \geqslant {I_{qh1}}\\ 
	&x{\kern 1pt} {\kern 1pt} {\kern 1pt} {\kern 1pt} {\kern 1pt} {\kern 1pt} {\kern 1pt} {\kern 1pt} {\kern 1pt} {\kern 1pt} {\kern 1pt} {\kern 1pt} {\kern 1pt} {\kern 1pt} {\kern 1pt} {\kern 1pt} {\kern 1pt} {\kern 1pt} {\kern 1pt} {\kern 1pt} {\kern 1pt} {\kern 1pt} {\kern 1pt} {\kern 1pt} {\kern 1pt} {\kern 1pt} {\kern 1pt} {\kern 1pt} {\kern 1pt} {\kern 1pt} {\kern 1pt} {\kern 1pt} {\kern 1pt} {\kern 1pt} {\kern 1pt} {\kern 1pt} {\kern 1pt} {\kern 1pt} if{\kern 1pt} {\kern 1pt} {\kern 1pt} {\kern 1pt} {I_{ql1}} \leqslant x \leqslant {I_{qh1}}\\ 
	&{I_{ql1 }}{\kern 1pt} {\kern 1pt} {\kern 1pt} {\kern 1pt} {\kern 1pt} {\kern 1pt} {\kern 1pt} {\kern 1pt} {\kern 1pt} {\kern 1pt} {\kern 1pt} {\kern 1pt} {\kern 1pt} {\kern 1pt} {\kern 1pt} {\kern 1pt} {\kern 1pt} {\kern 1pt} {\kern 1pt} {\kern 1pt} {\kern 1pt} {\kern 1pt} {\kern 1pt} {\kern 1pt} {\kern 1pt} {\kern 1pt} {\kern 1pt} {\kern 1pt} {\kern 1pt} {\kern 1pt} {\kern 1pt} {\kern 1pt} {\kern 1pt} {\kern 1pt} {\kern 1pt} {\kern 1pt} {\kern 1pt} {\kern 1pt} {\kern 1pt} {\kern 1pt} if{\kern 1pt} {\kern 1pt} {\kern 1pt} {\kern 1pt} x \leqslant {I_{ql1}}
	\end{aligned} \right.{\kern 1pt} {\kern 1pt} {\kern 1pt} {\kern 1pt} {\kern 1pt} {\kern 1pt} {\kern 1pt} {\kern 1pt} {\kern 1pt} {\kern 1pt}
	\end{eqnarray}
	\begin{eqnarray}\label{DBV}
	D{B_V}\left( {x} \right) = \left\{ \begin{aligned}
	& {x - dbd1} {\kern 1pt} {\kern 1pt} {\kern 1pt} {\kern 1pt} {\kern 1pt} {\kern 1pt} {\kern 1pt} {\kern 1pt} {\kern 1pt} {\kern 1pt} {\kern 1pt} {\kern 1pt} if{\kern 1pt} {\kern 1pt} {\kern 1pt} x> dbd1\\ 
	&0 {\kern 20pt} if{\kern 1pt} {\kern 1pt} {\kern 1pt} dbd2 \leqslant x \leqslant dbd1\\ 
	&{\kern 1pt}  {x - dbd2} {\kern 1pt} {\kern 1pt} {\kern 1pt} {\kern 1pt} {\kern 1pt} {\kern 1pt} {\kern 1pt} {\kern 1pt} if{\kern 1pt} {\kern 1pt} {\kern 1pt} x < dbd2{\kern 1pt} {\kern 1pt} {\kern 1pt} {\kern 1pt} {\kern 1pt} {\kern 1pt} {\kern 1pt} {\kern 1pt}
	\end{aligned} \right.
	\end{eqnarray}
	The current limit is modeled as follows:
	\begin{enumerate}
		\item Q-priority: $I_{qmax}=I_{max}$; $I_{qmin}=-I_{max}$;
		\item P-priority: $I_{qmax}=\sqrt{I_{max}^2-I_{qcmd}^2}$; if $typeflag=1$ then $I_{qmin}=-I_{qmax}$, else $I_{qmin}=0$.
	\end{enumerate}
	\subsubsection{Mathematical model of $S_4$}
	The local block diagram of $S_4$ is shown in Figure \ref{5}. The first block is a function of voltage tripping logic. Denoting it by a piecewise function as Equation (\ref{VT}), we can obtain the following dynamic equation:
	\begin{figure}[ht!]
		\centering
		\includegraphics[width=8cm]{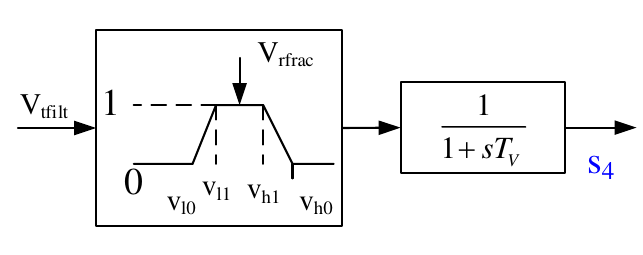}
		\caption{A local schematic of $S_4$ in the DER\_A model.}
		\label{5}
	\end{figure}
	\begin{figure}[ht!]
		\centering
		\includegraphics[width=7cm]{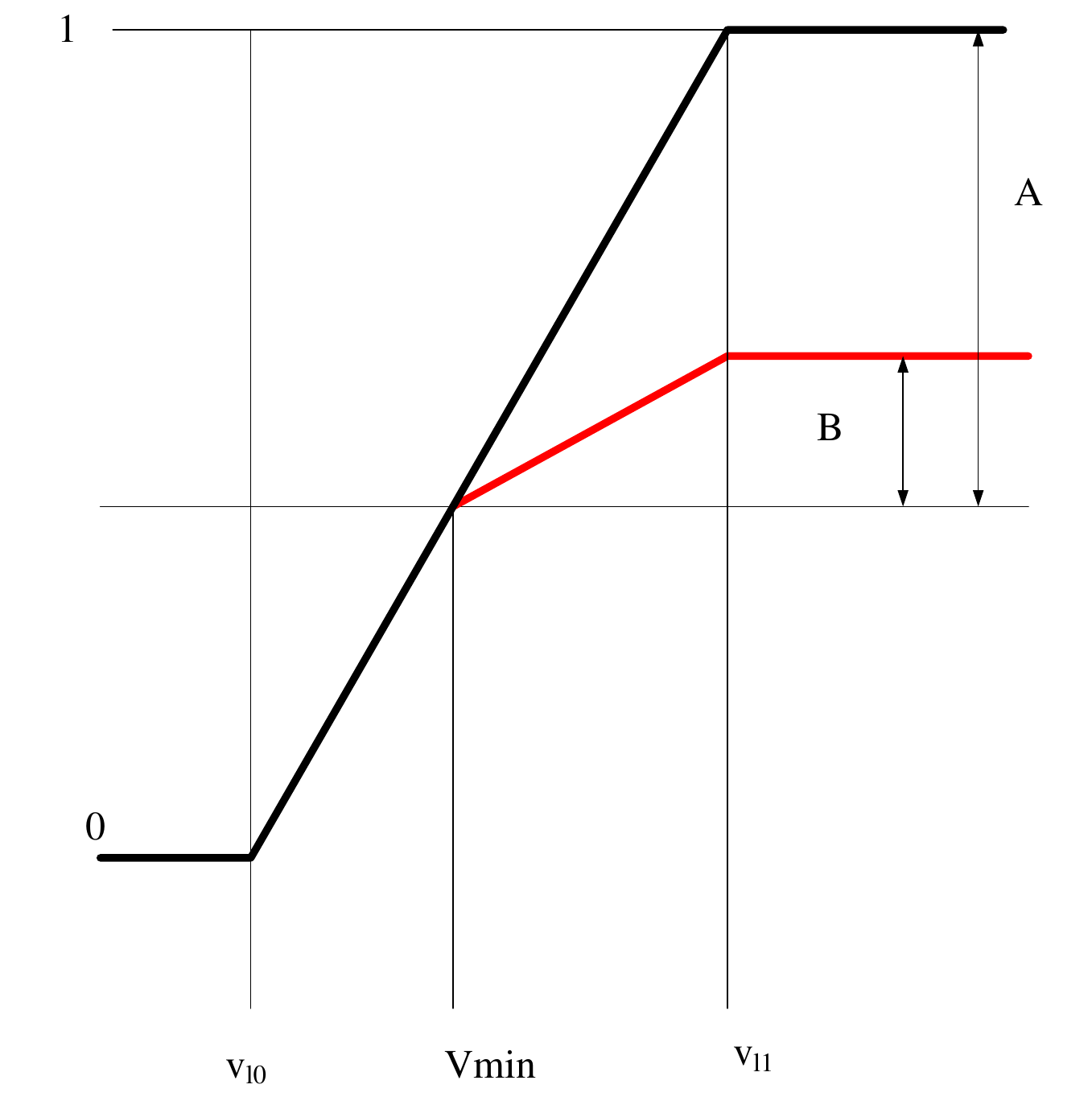}
		\caption{Effect of the Voltage trip}
		\label{voltagetrip}
	\end{figure}
	\begin{eqnarray}
	{\dot S_4} = \frac{1}{{{T_v}}}\left( {VoltageProtection({S_0},{V_{rfrac}}) - {S_4}} \right)
	\end{eqnarray}
	where the equations of voltage tripping logic is shown in Equation (\ref{VT}).
	\newcounter{mytempeqncnt}
	\begin{figure*}[!t]
		\normalsize
		\setcounter{mytempeqncnt}{\value{equation}}
		\setcounter{equation}{23}
		\begin{eqnarray}\label{VT}
		VoltageProtection({S_0},{V_{rfrac}}) = 
		\left\{ \begin{aligned}
		&\frac{{{V_t} - {V_{l0}}}}{{{V_{l1}} - {V_{l0}}}}{\kern 1pt} {\kern 1pt} {\kern 1pt} {\kern 1pt} {\kern 1pt} {\kern 1pt} {\kern 1pt} {\kern 1pt} {\kern 1pt} {\kern 1pt} {\kern 1pt} {\kern 1pt} {\kern 1pt} {\kern 1pt} {\kern 1pt} {\kern 1pt} {\kern 1pt} {\kern 1pt} {\kern 1pt} {\kern 1pt} {\kern 1pt} {\kern 1pt} {\kern 1pt} {\kern 1pt} {\kern 1pt} {\rm{ if}}{\kern 1pt} {\kern 1pt} {\kern 1pt} {\kern 1pt} {V_{l0}} \leqslant {{\rm{V}}_t} \leqslant {{\rm{V}}_{\min }}\\ 
		&\frac{{{V_t} - {V_{l0}}}}{{{V_{l1}} - {V_{l0}}}}{\kern 1pt} {\kern 1pt} {\kern 1pt} {\kern 1pt} {\kern 1pt} {\kern 1pt} {\kern 1pt} {\kern 1pt} {\kern 1pt} {\kern 1pt} {\kern 1pt} {\kern 1pt} {\kern 1pt} {\kern 1pt} {\kern 1pt} {\kern 1pt} {\kern 1pt} {\kern 1pt} {\kern 1pt} {\kern 1pt} {\kern 1pt} {\kern 1pt} {\kern 1pt} {\kern 1pt} {\kern 1pt} {\rm{ }}{\kern 1pt} {\rm{if}}{\kern 1pt} {\kern 1pt} {\kern 1pt} {\kern 1pt} {V_{\min }} \leqslant {{\rm{V}}_t} \leqslant {{\rm{V}}_{l1}}{\kern 1pt} {\kern 1pt} {\kern 1pt} {\kern 1pt} {\rm{and}}{\kern 1pt} {\kern 1pt} {\kern 1pt}  {{\rm{V}}_t} \leqslant {{\rm{V}}_{l1}} {\kern 1pt}{\kern 1pt}{\kern 1pt}\rm{for} {\kern 1pt}{\kern 1pt}{\kern 1pt}\rm{less}{\kern 1pt}{\kern 1pt}{\kern 1pt}\rm{than}{\kern 1pt}{\kern 1pt} {\kern 1pt}{{\rm{t}}_{lv1}}\\ 
		&1{\kern 1pt} {\kern 1pt} {\kern 1pt} {\kern 1pt} {\kern 1pt} {\kern 1pt} {\kern 1pt} {\kern 1pt} {\kern 1pt} {\kern 1pt} {\kern 1pt} {\kern 1pt} {\kern 1pt} {\kern 1pt} {\kern 1pt} {\kern 1pt} {\kern 1pt} {\kern 1pt} {\kern 1pt} {\kern 1pt} {\kern 1pt} {\kern 1pt} {\rm{ }}{\kern 1pt} {\kern 1pt} {\kern 1pt} {\kern 1pt} {\kern 1pt} {\kern 1pt} {\kern 1pt} {\kern 1pt} {\kern 1pt} {\kern 1pt} {\kern 1pt} {\kern 1pt} {\kern 1pt} {\kern 1pt} {\kern 1pt} {\kern 1pt} {\kern 1pt} {\kern 1pt} {\kern 1pt} {\kern 1pt} {\kern 1pt} {\kern 1pt} {\kern 1pt} {\kern 1pt} {\kern 1pt} {\kern 1pt} {\kern 1pt} {\kern 1pt} {\kern 1pt} {\kern 1pt} {\kern 1pt} {\kern 1pt} {\kern 1pt} {\kern 1pt} {\kern 1pt} {\kern 1pt} {\kern 1pt} {\kern 1pt} {\kern 1pt} {\kern 1pt}{\kern 1pt} {\kern 1pt}  {\rm{if}}{\kern 1pt} {\kern 1pt} {\kern 1pt} {\kern 1pt} {{\rm{V}}_{l1}} < {{\rm{V}}_t}{\rm{ < }}{\kern 1pt} {\kern 1pt} {V_{h1}}{\kern 1pt} {\kern 1pt} {\rm{and}}{\kern 1pt} {\kern 1pt} {\kern 1pt}  {{\rm{V}}_t} \leqslant {{\rm{V}}_{h1}} {\kern 1pt}{\kern 1pt}{\kern 1pt}\rm{for} {\kern 1pt}{\kern 1pt}{\kern 1pt}\rm{less}{\kern 1pt}{\kern 1pt}{\kern 1pt}\rm{than}{\kern 1pt}{\kern 1pt} {\kern 1pt}{{\rm{t}}_{lv1}}\\ 
		&\frac{{{V_{h0}} - {V_t}}}{{{V_{h0}} - {V_{h1}}}}{\kern 1pt} {\kern 1pt} {\kern 1pt} {\kern 1pt} {\kern 1pt} {\kern 1pt} {\kern 1pt} {\kern 1pt} {\kern 1pt} {\kern 1pt} {\kern 1pt} {\kern 1pt} {\kern 1pt} {\kern 1pt} {\kern 1pt} {\kern 1pt} {\kern 1pt} {\kern 1pt} {\kern 1pt} {\kern 1pt} {\kern 1pt} {\kern 1pt} {\kern 1pt} {\kern 1pt} {\kern 1pt} {\kern 1pt} {\kern 1pt} {\rm{if}}{\kern 1pt} {\kern 1pt} {\kern 1pt} {\kern 1pt} {V_{h1}} \leqslant {{\rm{V}}_t} \leqslant {V_{h0}}{\kern 1pt} {\kern 1pt} {\kern 1pt} {\kern 1pt} {\kern 1pt} {\rm{and}}{\kern 1pt} {\kern 1pt} {\kern 1pt}  {{\rm{V}}_t} \geqslant {{\rm{V}}_{h1}} {\kern 1pt}{\kern 1pt}{\kern 1pt}\rm{for} {\kern 1pt}{\kern 1pt}{\kern 1pt}\rm{less}{\kern 1pt}{\kern 1pt}{\kern 1pt}\rm{than}{\kern 1pt}{\kern 1pt} {\kern 1pt}{{\rm{t}}_{hv1}} \\ 
		&{V_{rfrac}}\frac{{{V_t} - {V_{\min }}}}{{{V_{l1}} - {V_{l0}}}}{\kern 1pt} {\kern 1pt} {\kern 1pt} {\kern 1pt} {\kern 1pt} {\kern 1pt} {\kern 1pt} {\rm{if}}{\kern 1pt} {\kern 1pt} {\kern 1pt} {\kern 1pt} {V_{\min }} \leqslant {{\rm{V}}_t} \leqslant {{\rm{V}}_{l1}}{\kern 1pt} {\kern 1pt} {\kern 1pt} {\kern 1pt} {\rm{and}}{\kern 1pt} {\kern 1pt} {\kern 1pt}  {{\rm{V}}_t} \leqslant {{\rm{V}}_{l1}} {\kern 1pt}{\kern 1pt}{\kern 1pt}\rm{for} {\kern 1pt}{\kern 1pt}{\kern 1pt}\rm{greater}{\kern 1pt}{\kern 1pt}{\kern 1pt}\rm{than}{\kern 1pt}{\kern 1pt} {\kern 1pt}{{\rm{t}}_{lv1}}\\ 
		&{V_{rfrac}}\left( {\frac{{{V_{l1}} - {V_{\min }}}}{{{V_{l1}} - {v_{l0}}}}} \right){\kern 1pt} {\kern 1pt} {\kern 1pt} {\kern 1pt} {\kern 1pt} {\kern 1pt} {\rm{if}}{\kern 1pt} {\kern 1pt} {\kern 1pt} {{\rm{V}}_{l1}} < {{\rm{V}}_t}{\rm{ < }}{\kern 1pt} {\kern 1pt} {V_{h1}}{\kern 1pt} {\kern 1pt} {\kern 1pt} {\kern 1pt} {\rm{and}}{\kern 1pt} {\kern 1pt} {\kern 1pt}  {{\rm{V}}_t} \leqslant {{\rm{V}}_{h1}} {\kern 1pt}{\kern 1pt}{\kern 1pt}\rm{for} {\kern 1pt}{\kern 1pt}{\kern 1pt}\rm{greater}{\kern 1pt}{\kern 1pt}{\kern 1pt}\rm{than}{\kern 1pt}{\kern 1pt} {\kern 1pt}{{\rm{t}}_{lv1}}\\ 
		&{V_{rfrac}}\left( {\frac{{{V_{\max }} - {V_t}}}{{{V_{h0}} - {V_{h1}}}}} \right){\kern 1pt} {\kern 1pt} {\kern 1pt} {\kern 1pt} {\kern 1pt} {\kern 1pt} {\rm{if}}{\kern 1pt} {\kern 1pt} {\kern 1pt} {{\rm{V}}_{h1}} \leqslant {{\rm{V}}_t} \leqslant {\kern 1pt} {\kern 1pt} {V_{\max }}{\kern 1pt} {\kern 1pt} {\kern 1pt} {\kern 1pt} {\rm{and}}{\kern 1pt} {\kern 1pt} {\kern 1pt}  {{\rm{V}}_t} \geqslant {{\rm{V}}_{h1}} {\kern 1pt}{\kern 1pt}{\kern 1pt}\rm{for} {\kern 1pt}{\kern 1pt}{\kern 1pt}\rm{greater}{\kern 1pt}{\kern 1pt}{\kern 1pt}\rm{than}{\kern 1pt}{\kern 1pt} {\kern 1pt}{{\rm{t}}_{hv1}}\\ 
		&\frac{{{V_{h0}} - {V_t}}}{{{V_{h0}} - {V_{h1}}}}{\kern 1pt} {\kern 1pt} {\kern 1pt} {\kern 1pt} {\kern 1pt} {\kern 1pt} {\kern 1pt} {\kern 1pt} {\kern 1pt} {\kern 1pt} {\kern 1pt} {\kern 1pt} {\kern 1pt} {\kern 1pt} {\kern 1pt} {\kern 1pt} {\kern 1pt} {\kern 1pt} {\kern 1pt} {\kern 1pt} {\kern 1pt} {\kern 1pt} {\kern 1pt} {\kern 1pt} {\kern 1pt} {\rm{ }}{\kern 1pt} {\kern 1pt} {\kern 1pt} {\kern 1pt} {\kern 1pt} {\kern 1pt} {\kern 1pt} {\kern 1pt} {\kern 1pt} {\kern 1pt} {\kern 1pt} {\kern 1pt} {\rm{if}}{\kern 1pt} {\kern 1pt} {\kern 1pt} {\kern 1pt} {V_{\max }} \leqslant {{\rm{V}}_t} \leqslant {{\rm{V}}_{h0}}\\ 
		&0{\kern 1pt} {\kern 1pt} {\kern 1pt} {\kern 1pt} {\kern 1pt} {\kern 1pt} {\kern 1pt} {\kern 1pt} {\kern 1pt} {\kern 1pt} {\kern 1pt} {\kern 1pt} {\kern 1pt} {\kern 1pt} {\kern 1pt} {\kern 1pt} {\kern 1pt} {\kern 1pt} {\kern 1pt} {\kern 1pt} {\kern 1pt} {\kern 1pt} {\kern 1pt} {\kern 1pt} {\kern 1pt} {\kern 1pt} {\kern 1pt} {\kern 1pt} {\kern 1pt} {\kern 1pt} {\kern 1pt} {\kern 1pt} {\kern 1pt} {\kern 1pt} {\kern 1pt} {\kern 1pt} {\kern 1pt} {\kern 1pt} {\kern 1pt} {\kern 1pt} {\kern 1pt} {\kern 1pt} {\kern 1pt} {\kern 1pt} {\kern 1pt} {\kern 1pt} {\kern 1pt} {\kern 1pt} {\kern 1pt} {\kern 1pt} {\kern 1pt} {\kern 1pt} {\kern 1pt} {\kern 1pt} {\kern 1pt} {\kern 1pt} {\kern 1pt} {\kern 1pt} {\kern 1pt} {\kern 1pt} {\kern 1pt} {\kern 1pt} {\kern 1pt} {\kern 1pt} {\kern 1pt} {\kern 1pt} {\kern 1pt} {\kern 1pt} {\kern 1pt} {\kern 1pt} {\kern 1pt} {\kern 1pt} {\kern 1pt} {\kern 1pt} {\kern 1pt} {\kern 1pt} otherwise
		\end{aligned} \right.
		\end{eqnarray}
		\setcounter{equation}{\value{mytempeqncnt}}
		\hrulefill
		\vspace*{4pt}
	\end{figure*}
	Note that $V_{min}$ is determined by internal software which keeps tracking
	the minimum voltage of $V_t$ during a simulation. Moreover, the frequency tripping logic is designed as follows: if frequency goes below $f_l$ for more than $t_{fl}$ seconds, then the entire model will trip; if frequency goes above $f_h$ for more than $t_{fh}$ seconds, then the entire model will trip.
	\subsubsection{Mathematical model of $S_5$}
	\begin{figure}[ht!]
		\centering
		\includegraphics[width=7cm]{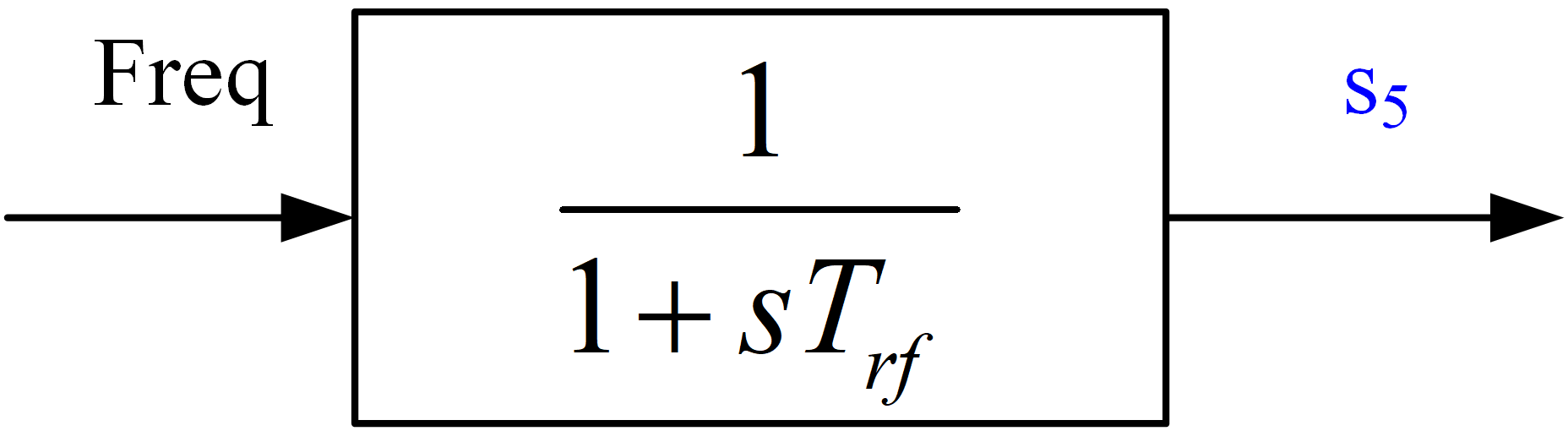}
		\caption{A local schematic of $S_5$ in the DER\_A model.}
		\label{6}
	\end{figure}
	Figure \ref{6} shows the block diagram of first-order filter whose input is the input frequency $Freq$, and the output is filtered frequency $S_5$ ($Freq_{filt}$). From the diagram, we can obtain the following dynamic equation:
	\setcounter{equation}{24}
	\begin{eqnarray}
	{\dot S_5} = \frac{1}{{{T_{rf}}}}\left( {Freq - {S_5}} \right)
	\end{eqnarray}
	\subsubsection{Mathematical model of $S_6$}
	Figure \ref{7} shows the diagram of PI controller with respect to $S_6$. Defining the limiter and dead bands functions as Equation (\ref{sat3}) - (\ref{Gup}), we can obtain the following model of $S_6$:
	\begin{figure}[ht!]
		\centering
		\includegraphics[width=8cm]{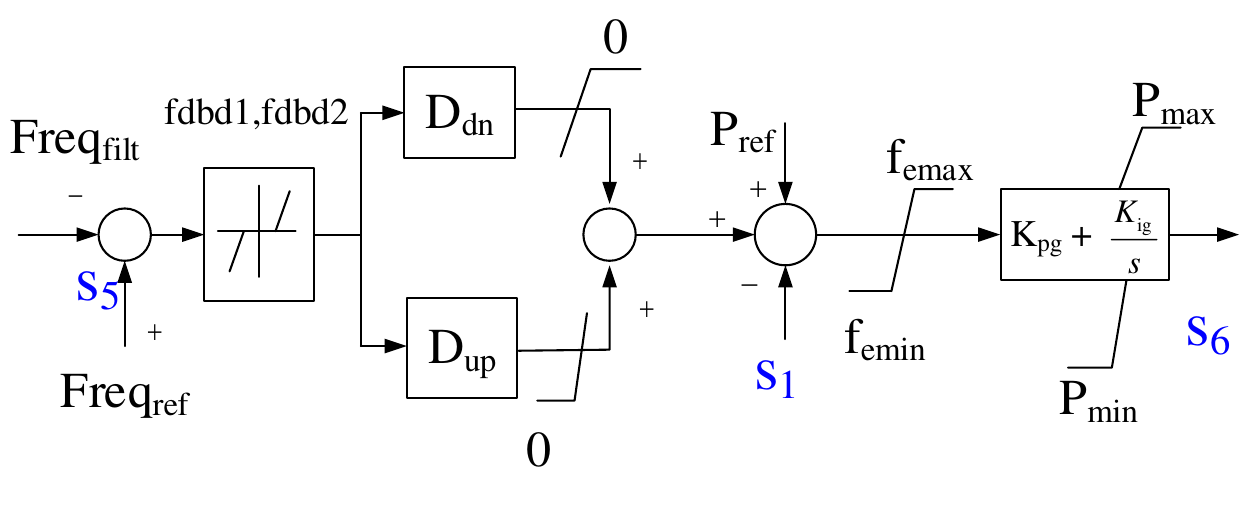}
		\caption{A local schematic of $S_6$ in the DER\_A model.}
		\label{7}
	\end{figure}
	\begin{equation}\label{s6}
	\begin{aligned}
	{\dot S_6} \!= &{K_{ig}}sa{t_4}( {P_{ref}} \!- \!{S_1}\! +\! sa{t_5}\left[ {{D_{dn}} \!\cdot \!D{B_F}(Fre{q_{ref}}\! -\! {S_5})} \right] \!\!\!\!\!\!\!\!\!\!\!\!\!\!\!\\\!\!\!\!\!\!\!\!
	&+\! sa{t_6}\left[ {{D_{up}} \!\cdot\! D{B_F}(Fre{q_{ref}}\! -\! {S_5})} \right] )\!+\! {\frac{K_{pg}}{T_p}}{S_1}\!\!\!\!\!\!\!\!\!\!\!\!\!\!\!\\\!\!\!\!\!\!\!\!
	&+\!G_{dn}\left( {Freq-S_5}\right)+G_{up}\left({Freq-S_5}\right)  -\frac{S_{8}}{T_p}\!\!\!\!\!\!\!\!
	\end{aligned}
	\end{equation}
	\begin{equation}
	sa{t_4}\left( x \right) = \left\{ \begin{aligned}\label{sat4}
	&{f_{e\max }}{\kern 1pt} {\kern 1pt} {\kern 1pt} {\kern 1pt} {\kern 1pt} {\kern 1pt} {\kern 1pt} {\kern 1pt} {\kern 1pt} {\kern 1pt} {\kern 1pt} {\kern 1pt} {\kern 1pt} {\kern 1pt} {\kern 1pt} {\kern 1pt} {\kern 1pt} {\kern 1pt} {\kern 1pt} {\kern 1pt} {\kern 1pt} {\kern 1pt} {\kern 1pt} {\kern 1pt} {\kern 1pt} {\kern 1pt} {\kern 1pt} {\kern 1pt} {\kern 1pt} {\kern 1pt} {\kern 1pt} {\kern 1pt} {\kern 1pt} {\kern 1pt} {\kern 1pt} {\kern 1pt} {\kern 1pt} {\kern 1pt} if{\kern 1pt} {\kern 1pt} {\kern 1pt} x \geqslant {f_{e\max }}\\ 
	&x{\kern 1pt} {\kern 1pt} {\kern 1pt} {\kern 1pt} {\kern 1pt} {\kern 1pt} {\kern 1pt} {\kern 1pt} {\kern 1pt} {\kern 1pt} {\kern 1pt} {\kern 1pt} {\kern 1pt} {\kern 1pt} {\kern 1pt} {\kern 1pt} {\kern 1pt} {\kern 1pt} {\kern 1pt} {\kern 1pt} {\kern 1pt} {\kern 1pt} {\kern 1pt} {\kern 1pt} {\kern 1pt} {\kern 1pt} {\kern 1pt} {\kern 1pt} {\kern 1pt} {\kern 1pt} {\kern 1pt} {\kern 1pt} {\kern 1pt} {\kern 1pt} {\kern 1pt} {\kern 1pt} {\kern 1pt} {\kern 1pt} if{\kern 1pt} {\kern 1pt} {\kern 1pt} {\kern 1pt} {f_{e\min }} \leqslant x \leqslant {f_{e\max }}\\ 
	&{f_{e\min }}{\kern 1pt} {\kern 1pt} {\kern 1pt} {\kern 1pt} {\kern 1pt} {\kern 1pt} {\kern 1pt} {\kern 1pt} {\kern 1pt} {\kern 1pt} {\kern 1pt} {\kern 1pt} {\kern 1pt} {\kern 1pt} {\kern 1pt} {\kern 1pt} {\kern 1pt} {\kern 1pt} {\kern 1pt} {\kern 1pt} {\kern 1pt} {\kern 1pt} {\kern 1pt} {\kern 1pt} {\kern 1pt} {\kern 1pt} {\kern 1pt} {\kern 1pt} {\kern 1pt} {\kern 1pt} {\kern 1pt} {\kern 1pt} {\kern 1pt} {\kern 1pt} {\kern 1pt} {\kern 1pt} {\kern 1pt} {\kern 1pt} {\kern 1pt} {\kern 1pt} if{\kern 1pt} {\kern 1pt} {\kern 1pt} {\kern 1pt} x \leqslant {f_{e\min }}
	\end{aligned} \right.{\kern 1pt} {\kern 1pt} {\kern 1pt} {\kern 1pt} {\kern 1pt} {\kern 1pt} {\kern 1pt} {\kern 1pt} {\kern 1pt} {\kern 1pt}
	\end{equation}
	\begin{equation}\label{sat5}
	sa{t_5}\left( x \right) = \left\{ \begin{aligned}
	&x{\kern 1pt} {\kern 1pt} {\kern 1pt} {\kern 1pt} {\kern 1pt} {\kern 1pt} {\kern 1pt} {\kern 1pt} {\kern 1pt} {\kern 1pt} {\kern 1pt} {\kern 1pt} {\kern 1pt} {\kern 1pt} {\kern 1pt} {\kern 1pt} {\kern 1pt} {\kern 1pt} {\kern 1pt} {\kern 1pt} {\kern 1pt} {\kern 1pt} {\kern 1pt} {\kern 1pt} {\kern 1pt} {\kern 1pt} {\kern 1pt} {\kern 1pt} {\kern 1pt} {\kern 1pt} if{\kern 1pt} {\kern 1pt} {\kern 1pt} x \leqslant 0\\ 
	&0{\kern 1pt} {\kern 1pt} {\kern 1pt} {\kern 1pt} {\kern 1pt} {\kern 1pt} {\kern 1pt} {\kern 1pt} {\kern 1pt} {\kern 1pt} {\kern 1pt} {\kern 1pt} {\kern 1pt} {\kern 1pt} {\kern 1pt} {\kern 1pt} {\kern 1pt} {\kern 1pt} {\kern 1pt} {\kern 1pt} {\kern 1pt} {\kern 1pt} {\kern 1pt} {\kern 1pt} {\kern 1pt} {\kern 1pt} {\kern 1pt} {\kern 1pt} {\kern 1pt} {\kern 1pt} if{\kern 1pt} {\kern 1pt} {\kern 1pt} {\kern 1pt} x > 0
	\end{aligned} \right.{\kern 1pt} {\kern 1pt} {\kern 1pt} {\kern 1pt} {\kern 1pt} {\kern 1pt} {\kern 1pt} {\kern 1pt} {\kern 1pt} {\kern 1pt}
	\end{equation}
	\begin{equation}\label{sat6}
	sa{t_6}\left( x \right) = \left\{ \begin{aligned}
	&x{\kern 1pt} {\kern 1pt} {\kern 1pt} {\kern 1pt} {\kern 1pt} {\kern 1pt} {\kern 1pt} {\kern 1pt} {\kern 1pt} {\kern 1pt} {\kern 1pt} {\kern 1pt} {\kern 1pt} {\kern 1pt} {\kern 1pt} {\kern 1pt} {\kern 1pt} {\kern 1pt} {\kern 1pt} {\kern 1pt} {\kern 1pt} {\kern 1pt} {\kern 1pt} {\kern 1pt} {\kern 1pt} {\kern 1pt} {\kern 1pt} {\kern 1pt} {\kern 1pt} {\kern 1pt} if{\kern 1pt} {\kern 1pt} {\kern 1pt} x > 0\\ 
	&0{\kern 1pt} {\kern 1pt} {\kern 1pt} {\kern 1pt} {\kern 1pt} {\kern 1pt} {\kern 1pt} {\kern 1pt} {\kern 1pt} {\kern 1pt} {\kern 1pt} {\kern 1pt} {\kern 1pt} {\kern 1pt} {\kern 1pt} {\kern 1pt} {\kern 1pt} {\kern 1pt} {\kern 1pt} {\kern 1pt} {\kern 1pt} {\kern 1pt} {\kern 1pt} {\kern 1pt} {\kern 1pt} {\kern 1pt} {\kern 1pt} {\kern 1pt} {\kern 1pt} {\kern 1pt} if{\kern 1pt} {\kern 1pt} {\kern 1pt} {\kern 1pt} x \leqslant 0
	\end{aligned} \right.{\kern 1pt} {\kern 1pt} {\kern 1pt} {\kern 1pt} {\kern 1pt} {\kern 1pt} {\kern 1pt} {\kern 1pt} {\kern 1pt} {\kern 1pt}
	\end{equation}
	\begin{equation}\label{DBF}
	D{B_F}\left( {x} \right) = \left\{ \begin{aligned}
	& {x - fdbd2} {\kern 1pt} {\kern 1pt} {\kern 1pt} {\kern 1pt} {\kern 1pt} {\kern 1pt} {\kern 1pt} {\kern 1pt} {\kern 1pt} {\kern 1pt} {\kern 1pt} {\kern 1pt} if{\kern 1pt} {\kern 1pt} {\kern 1pt} x> fdbd2\\ 
	&0{\kern 1pt} {\kern 1pt} {\kern 1pt} {\kern 1pt} {\kern 1pt} {\kern 1pt} {\kern 1pt} {\kern 1pt} {\kern 1pt} {\kern 1pt} {\kern 1pt} {\kern 1pt} {\kern 1pt} {\kern 1pt} {\kern 1pt} {\kern 1pt} {\kern 1pt} {\kern 1pt} {\kern 1pt} {\kern 1pt} {\kern 1pt} {\kern 1pt} {\kern 1pt} {\kern 1pt} {\kern 1pt} {\kern 1pt} {\kern 1pt} {\kern 1pt} {\kern 1pt} {\kern 1pt} {\kern 1pt} {\kern 1pt} {\kern 1pt} {\kern 1pt} {\kern 1pt} {\kern 1pt} {\kern 1pt} {\kern 1pt} {\kern 1pt} {\kern 1pt} {\kern 1pt} {\kern 1pt} {\kern 1pt}  if{\kern 1pt} {\kern 1pt} {\kern 1pt} fdbd1 \leqslant x \leqslant fdbd2\\ 
	&{\kern 1pt}  {x - fdbd1} {\kern 1pt} {\kern 1pt} {\kern 1pt} {\kern 1pt} {\kern 1pt} {\kern 1pt} {\kern 1pt} {\kern 1pt} if{\kern 1pt} {\kern 1pt} {\kern 1pt} x < fdbd1{\kern 1pt} {\kern 1pt} {\kern 1pt} {\kern 1pt} {\kern 1pt} {\kern 1pt} {\kern 1pt} {\kern 1pt}
	\end{aligned} \right.
	\end{equation}
	\begin{equation}\label{Gdn}
	G_{dn}\left( {x}\right) =\left\{ \begin{aligned}
	&{-\frac{{K_{pg}}{D_{dn}}}{T_{rf}}}{x} {\kern 10pt} {{if}}{\kern 1pt} {\kern 1pt} {\kern 1pt}{x}{<} fdbd1{\kern 1pt} {\kern 1pt} {\kern 1pt}{\kern 1pt} {\kern 1pt} {\kern 1pt}{or}{\kern 1pt} {\kern 1pt} {\kern 1pt}{x}{>}fdb2, {\kern 1pt} {\kern 1pt} {\kern 1pt}\!\!\!\!\!\!\!\!\!\\\!\!\!\!\!\!\!\!\!
	&{\kern 65pt}{and}{\kern 1pt} {\kern 1pt} {\kern 1pt}{\frac{{D_{dn}}}{T_{rf}}} {x} \geqslant 0\!\!\!\!\!\!\!\!\!\\\!\!\!\!\!\!\!\!\!
	&0 {\kern 1pt} {\kern 1pt} {\kern 1pt} {\kern 1pt} {\kern 1pt}{\kern 1pt} {\kern 1pt} {\kern 1pt} {\kern 1pt} {\kern 1pt} {\kern 1pt} {\kern 1pt} {\kern 1pt} {\kern 1pt} {\kern 1pt}{\kern 1pt} {\kern 1pt} {\kern 1pt} {\kern 1pt} {\kern 1pt} {\kern 1pt} {\kern 1pt} {\kern 1pt} {\kern 1pt} {\kern 1pt}{\kern 1pt} {\kern 1pt} {\kern 1pt} {\kern 1pt} {\kern 1pt} {\kern 1pt} {\kern 1pt} {\kern 1pt} {\kern 1pt} {\kern 1pt}{\kern 1pt} {\kern 1pt} {\kern 1pt} {\kern 1pt} {\kern 1pt} {\kern 1pt} {\kern 1pt} {\kern 1pt} {\kern 1pt} {\kern 1pt}{\kern 1pt} {\kern 1pt} {\kern 1pt} {\kern 1pt} {\kern 1pt} {\kern 1pt} {\kern 1pt} {\kern 1pt} {\kern 1pt} {\kern 1pt}{\kern 1pt} {\kern 1pt} {\kern 1pt}{\kern 1pt} {\kern 1pt} {\kern 1pt}\rm{otherwise}\!\!\!\!\!\!\!\!\!
	\end{aligned} \right.
	\end{equation}
	\begin{equation}\label{Gup}
	G_{up}\left( {x}\right) =\left\{ \begin{aligned}
	&{-\frac{{K_{pg}}{D_{up}}}{T_{rf}}}{x}{\kern 10pt} {{if}}{\kern 1pt} {\kern 1pt} {\kern 1pt}{x}{<} fdbd1{\kern 1pt} {\kern 1pt} {\kern 1pt}{\kern 1pt} {\kern 1pt} {\kern 1pt}{or}{\kern 1pt} {\kern 1pt} {\kern 1pt}{x}{>}fdb2, {\kern 1pt} {\kern 1pt} {\kern 1pt}\!\!\!\!\!\!\!\!\!\\\!\!\!\!\!\!\!\!\!
	&{\kern 65pt}{and}{\kern 1pt} {\kern 1pt} {\kern 1pt}{\frac{{D_{up}}}{T_{rf}}} {x} < 0\!\!\!\!\!\!\!\!\!\\\!\!\!\!\!\!\!\!\!
	&0{\kern 1pt} {\kern 1pt} {\kern 1pt} {\kern 1pt} {\kern 1pt}{\kern 1pt} {\kern 1pt} {\kern 1pt} {\kern 1pt} {\kern 1pt} {\kern 1pt} {\kern 1pt} {\kern 1pt} {\kern 1pt} {\kern 1pt}{\kern 1pt} {\kern 1pt} {\kern 1pt} {\kern 1pt} {\kern 1pt} {\kern 1pt} {\kern 1pt} {\kern 1pt} {\kern 1pt} {\kern 1pt}{\kern 1pt} {\kern 1pt} {\kern 1pt} {\kern 1pt} {\kern 1pt} {\kern 1pt} {\kern 1pt} {\kern 1pt} {\kern 1pt} {\kern 1pt}{\kern 1pt} {\kern 1pt} {\kern 1pt} {\kern 1pt} {\kern 1pt} {\kern 1pt} {\kern 1pt} {\kern 1pt} {\kern 1pt} {\kern 1pt}{\kern 1pt} {\kern 1pt} {\kern 1pt} {\kern 1pt} {\kern 1pt} {\kern 1pt} {\kern 1pt} {\kern 1pt} {\kern 1pt} {\kern 1pt}{\kern 1pt} {\kern 1pt} {\kern 1pt}{\kern 1pt} {\kern 1pt} {\kern 1pt}\rm{otherwise}\!\!\!\!\!\!\!\!\!
	\end{aligned} \right.
	\end{equation}
	\subsubsection{Mathematical model of $S_7$}
	The local block diagram of $S_7$ is shown in Figure \ref{8}. From the diagram, we can obtain the following dynamic equation:
	\begin{figure}[ht!]
		\centering
		\includegraphics[width=7cm]{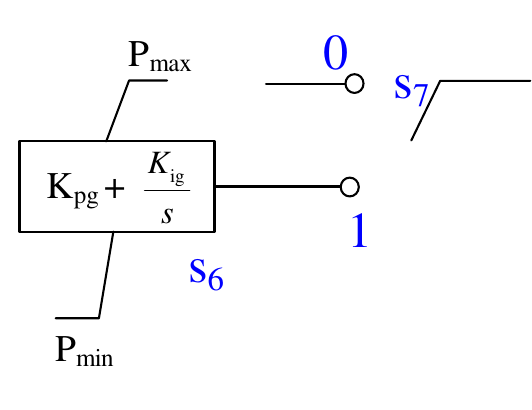}
		\caption{A local schematic of $S_{7}$ in the DER\_A model.}
		\label{8}
	\end{figure}
	\begin{equation}
	{\dot S_7} = \left\{ \begin{aligned}
	&{0}{\kern 68pt} {\kern 1pt} {\kern 1pt} {\kern 1pt} {\kern 1pt} {\kern 1pt} {\kern 1pt} {\kern 1pt} {\kern 1pt} {\kern 1pt} {\kern 1pt} {\kern 1pt} {\kern 1pt} {\kern 1pt} {\kern 1pt} {\kern 1pt} {\kern 1pt} {\kern 1pt} {\kern 1pt} if{\kern 1pt} {\kern 1pt} {\kern 1pt} Fre{q_{flag}} = 0\\
	&sat_8\left(s\dot{a}{t_7}\left( {{S_6}} \right)\right){\kern 1pt} {\kern 1pt} {\kern 1pt} {\kern 1pt} {\kern 1pt} {\kern 1pt} {\kern 1pt} {\kern 1pt} {\kern 1pt} {\kern 1pt} {\kern 1pt} {\kern 1pt} {\kern 1pt} {\kern 1pt} {\kern 1pt} {\kern 1pt} {\kern 1pt} {\kern 1pt} {\kern 1pt} {\kern 1pt} {\kern 1pt} {\kern 1pt} {\kern 1pt} if{\kern 1pt} {\kern 1pt} {\kern 1pt} {\kern 1pt} Fre{q_{flag}} = 1
	\end{aligned} \right.{\kern 1pt} {\kern 1pt} {\kern 1pt} {\kern 1pt} {\kern 1pt} {\kern 1pt} {\kern 1pt} {\kern 1pt} {\kern 1pt} {\kern 1pt}
	\end{equation}
	where the limiter function is defined as follows,
	\begin{equation}\label{sat7}
	sa{t_7}\left( x \right) = \left\{ \begin{aligned}
	&{P_{\max }}{\kern 1pt} {\kern 1pt} {\kern 1pt} {\kern 1pt} {\kern 1pt} {\kern 1pt} {\kern 1pt} {\kern 1pt} {\kern 1pt} {\kern 1pt} {\kern 1pt} {\kern 1pt} {\kern 1pt} {\kern 1pt} {\kern 1pt} {\kern 1pt} {\kern 1pt} {\kern 1pt} {\kern 1pt} {\kern 1pt} {\kern 1pt} {\kern 1pt} {\kern 1pt} {\kern 1pt} {\kern 1pt} {\kern 1pt} {\kern 1pt} {\kern 1pt} {\kern 1pt} {\kern 1pt} {\kern 1pt} {\kern 1pt} {\kern 1pt} {\kern 1pt} {\kern 1pt} {\kern 1pt} {\kern 1pt} {\kern 1pt} if{\kern 1pt} {\kern 1pt} {\kern 1pt} x \geqslant {P_{\max }}\\ 
	&x{\kern 1pt} {\kern 1pt} {\kern 1pt} {\kern 1pt} {\kern 1pt} {\kern 1pt} {\kern 1pt} {\kern 1pt} {\kern 1pt} {\kern 1pt} {\kern 1pt} {\kern 1pt} {\kern 1pt} {\kern 1pt} {\kern 1pt} {\kern 1pt} {\kern 1pt} {\kern 1pt} {\kern 1pt} {\kern 1pt} {\kern 1pt} {\kern 1pt} {\kern 1pt} {\kern 1pt} {\kern 1pt} {\kern 1pt} {\kern 1pt} {\kern 1pt} {\kern 1pt} {\kern 1pt} {\kern 1pt} {\kern 1pt} {\kern 1pt} {\kern 1pt} {\kern 1pt} {\kern 1pt} {\kern 1pt} {\kern 1pt} if{\kern 1pt} {\kern 1pt} {\kern 1pt} {\kern 1pt} {P_{\min }} \leqslant x \leqslant {P_{\max }}\\ 
	&{P_{\min }}{\kern 1pt} {\kern 1pt} {\kern 1pt} {\kern 1pt} {\kern 1pt} {\kern 1pt} {\kern 1pt} {\kern 1pt} {\kern 1pt} {\kern 1pt} {\kern 1pt} {\kern 1pt} {\kern 1pt} {\kern 1pt} {\kern 1pt} {\kern 1pt} {\kern 1pt} {\kern 1pt} {\kern 1pt} {\kern 1pt} {\kern 1pt} {\kern 1pt} {\kern 1pt} {\kern 1pt} {\kern 1pt} {\kern 1pt} {\kern 1pt} {\kern 1pt} {\kern 1pt} {\kern 1pt} {\kern 1pt} {\kern 1pt} {\kern 1pt} {\kern 1pt} {\kern 1pt} {\kern 1pt} {\kern 1pt} {\kern 1pt} {\kern 1pt} {\kern 1pt} if{\kern 1pt} {\kern 1pt} {\kern 1pt} {\kern 1pt} x \leqslant {P_{\min }}
	\end{aligned} \right.{\kern 1pt} {\kern 1pt} {\kern 1pt} {\kern 1pt} {\kern 1pt} {\kern 1pt} {\kern 1pt} {\kern 1pt} {\kern 1pt} {\kern 1pt}
	\end{equation}
	\begin{equation}\label{sat8}
	sa{t_8}\left( x \right) = \left\{ \begin{aligned}
	&d{P_{\max }}{\kern 1pt} {\kern 1pt} {\kern 1pt} {\kern 1pt} {\kern 1pt} {\kern 1pt} {\kern 1pt} {\kern 1pt} {\kern 1pt} {\kern 1pt} {\kern 1pt} {\kern 1pt} {\kern 1pt} {\kern 1pt} {\kern 1pt} {\kern 1pt} {\kern 1pt} {\kern 1pt} {\kern 1pt} {\kern 1pt} {\kern 1pt} {\kern 1pt} {\kern 1pt} {\kern 1pt} {\kern 1pt} {\kern 1pt} {\kern 1pt} {\kern 1pt} {\kern 1pt} {\kern 1pt} {\kern 1pt} {\kern 1pt} {\kern 1pt} {\kern 1pt} {\kern 1pt} {\kern 1pt} {\kern 1pt} {\kern 1pt} if{\kern 1pt} {\kern 1pt} {\kern 1pt} x \geqslant d{P_{\max }}\\ 
	&x{\kern 1pt} {\kern 1pt} {\kern 1pt} {\kern 1pt} {\kern 1pt} {\kern 1pt} {\kern 1pt} {\kern 1pt} {\kern 1pt} {\kern 1pt} {\kern 1pt} {\kern 1pt} {\kern 1pt} {\kern 1pt} {\kern 1pt} {\kern 1pt} {\kern 1pt} {\kern 1pt} {\kern 1pt} {\kern 1pt} {\kern 1pt} {\kern 1pt} {\kern 1pt} {\kern 1pt} {\kern 1pt} {\kern 1pt} {\kern 1pt} {\kern 1pt} {\kern 1pt} {\kern 1pt} {\kern 1pt} {\kern 1pt} {\kern 1pt} {\kern 1pt} {\kern 1pt} {\kern 1pt} {\kern 1pt} {\kern 1pt} if{\kern 1pt} {\kern 1pt} {\kern 1pt} {\kern 1pt} d{P_{\min }} \leqslant x \leqslant d{P_{\max }}\\ 
	&d{P_{\min }}{\kern 1pt} {\kern 1pt} {\kern 1pt} {\kern 1pt} {\kern 1pt} {\kern 1pt} {\kern 1pt} {\kern 1pt} {\kern 1pt} {\kern 1pt} {\kern 1pt} {\kern 1pt} {\kern 1pt} {\kern 1pt} {\kern 1pt} {\kern 1pt} {\kern 1pt} {\kern 1pt} {\kern 1pt} {\kern 1pt} {\kern 1pt} {\kern 1pt} {\kern 1pt} {\kern 1pt} {\kern 1pt} {\kern 1pt} {\kern 1pt} {\kern 1pt} {\kern 1pt} {\kern 1pt} {\kern 1pt} {\kern 1pt} {\kern 1pt} {\kern 1pt} {\kern 1pt} {\kern 1pt} {\kern 1pt} {\kern 1pt} {\kern 1pt} {\kern 1pt} if{\kern 1pt} {\kern 1pt} {\kern 1pt} {\kern 1pt} x \leqslant d{P_{\min }}
	\end{aligned} \right.{\kern 1pt} {\kern 1pt} {\kern 1pt} {\kern 1pt} {\kern 1pt} {\kern 1pt} {\kern 1pt} {\kern 1pt} {\kern 1pt} {\kern 1pt}
	\end{equation}
	\subsubsection{Mathematical model of $S_8$}
	The local block diagram of $S_8$ is shown in Figure \ref{9}. From the diagram, we can obtain the following dynamic equation:
	\begin{figure}[ht!]
		\centering
		\includegraphics[width=7cm]{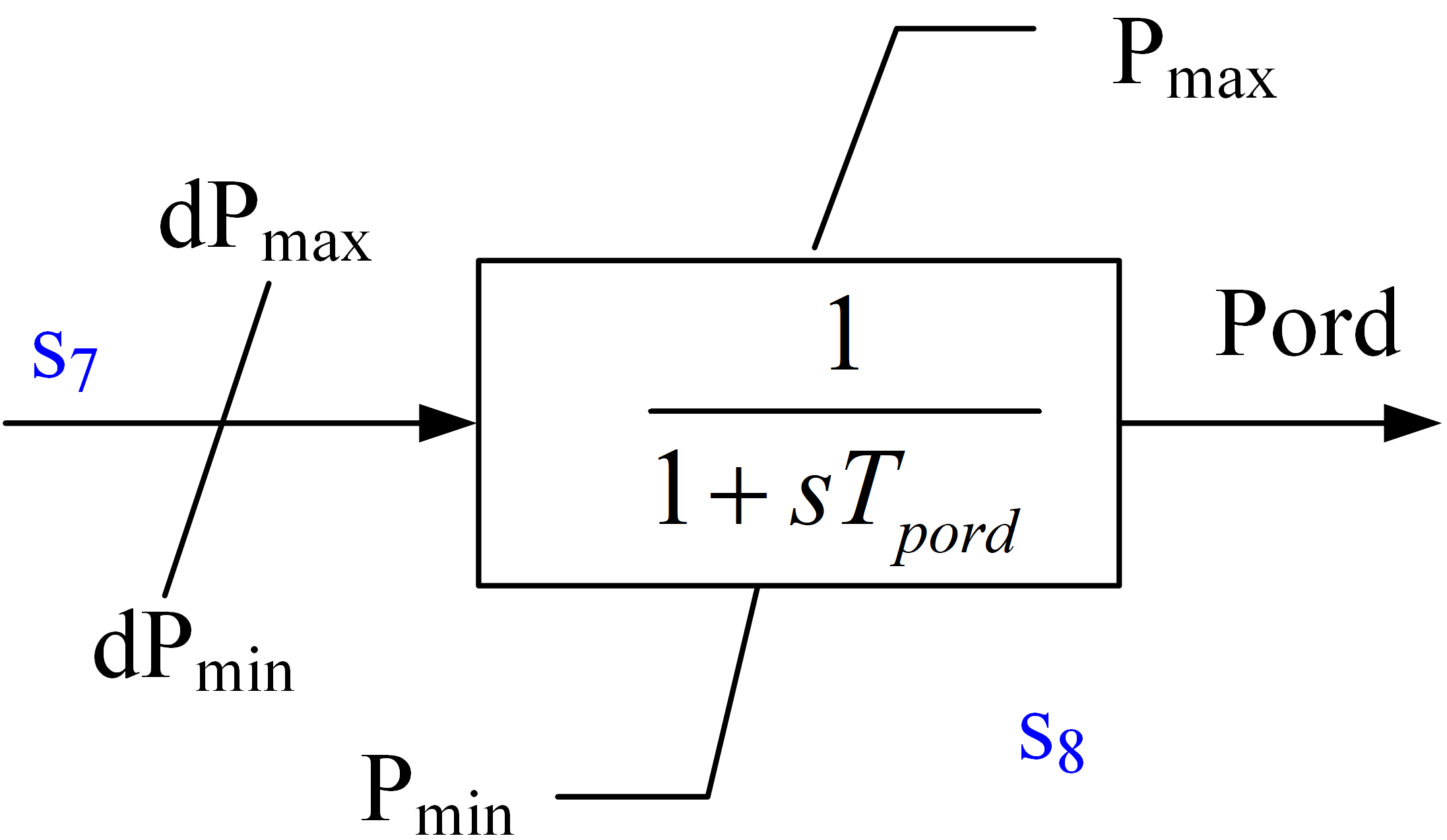}
		\caption{A local schematic of $S_{8}$ in the DER\_A model.}
		\label{9}
	\end{figure}
	\begin{equation}
	{\dot S_8} = \frac{1}{{{T_{pord}}}}\left( {{S_7} - {S_8}} \right)
	\end{equation}
	\subsubsection{Mathematical model of $S_9$}
	The local block diagram of d-axis current $S_9$ ($i_d$) is shown in Figure \ref{10}. From the diagram, we can obtain the following dynamic equation:
	\begin{figure}[ht!]
		\centering
		\includegraphics[width=8cm]{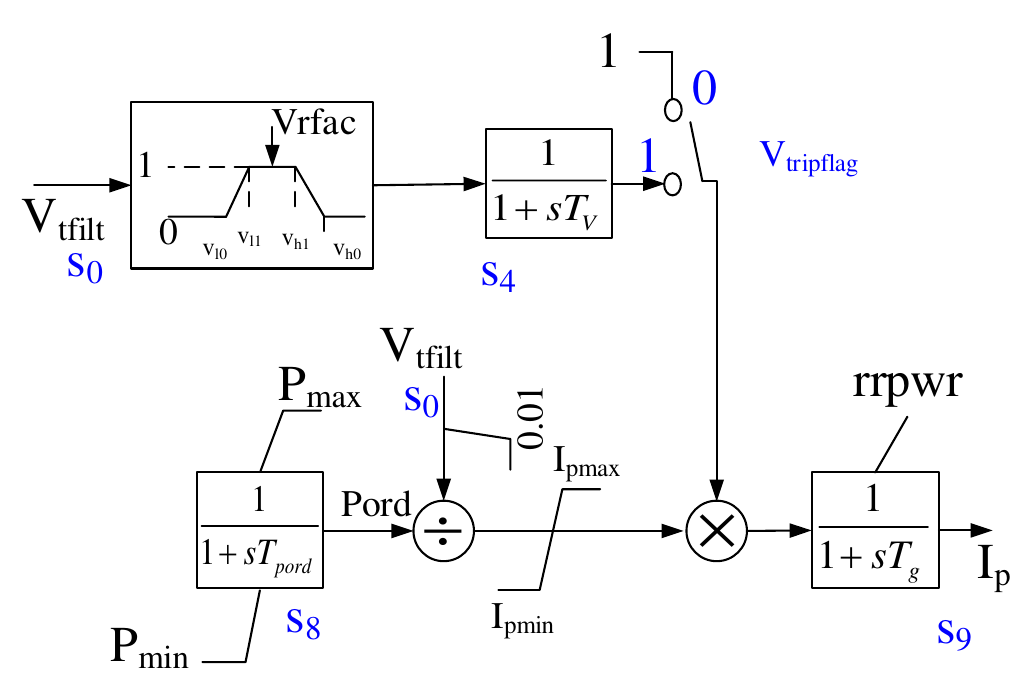}
		\caption{A local schematic of $S_{9}$ in the DER\_A model.}
		\label{10}
	\end{figure}
	\begin{eqnarray}
	{\dot S_9}\! =\! \left\{\! \begin{aligned}
	&\!\frac{1}{T_g}\!\left(\! {sat_{9}\!\left( {\frac{{sat_7({S_8})}}{{sat_1\left( S_0 \right)}}}\right) \!\!\times\!\! {S_4}\!\! -\!\! {S_9}} \right){\kern 2pt} if {\kern 2pt} {V_{tripflag}}\! =\! 1\!\!\!\!\!\\\!\!\!\!\!\!
	\!\!\!\!\!\\\!\!\!\!\!
	&\!\frac{1}{T_g}\!\left(\! {sat_{9}\!\left( {\frac{{sat_7({S_8})}}{{sa{t_1}\left( {{S_0}} \right)}}}\right) \!-\! {S_9}}\right){\kern 15pt} if  {\kern 5pt} {V_{tripflag}} \!=\! 0\!\!\!\!\!
	\end{aligned} \right.
	\end{eqnarray}
	where the limiter function is defined as Equation (\ref{sat1}), (\ref{sat7}) and (\ref{sat9}). 
	\begin{equation}\label{sat9}
	sa{t_9}\left( x \right) = \left\{ \begin{aligned}
	&{I_{p\max }} {\kern 1pt} {\kern 1pt} {\kern 1pt} {\kern 1pt} {\kern 1pt} {\kern 1pt} {\kern 1pt} {\kern 1pt} {\kern 1pt} {\kern 1pt} {\kern 1pt} {\kern 1pt} {\kern 1pt} {\kern 1pt} {\kern 1pt} {\kern 1pt} {\kern 1pt} {\kern 1pt} {\kern 1pt} {\kern 1pt} {\kern 1pt} {\kern 1pt} {\kern 1pt} {\kern 1pt} {\kern 1pt} {\kern 1pt} {\kern 1pt} {\kern 1pt} {\kern 1pt} {\kern 1pt} {\kern 1pt} {\kern 1pt} {\kern 1pt} if{\kern 1pt} {\kern 1pt} {\kern 1pt} x \geqslant {I_{p\max }}\\ 
	&x{\kern 1pt} {\kern 1pt} {\kern 1pt} {\kern 1pt} {\kern 1pt} {\kern 1pt} {\kern 1pt} {\kern 1pt} {\kern 1pt} {\kern 1pt} {\kern 1pt} {\kern 1pt} {\kern 1pt} {\kern 1pt} {\kern 1pt} {\kern 1pt} {\kern 1pt} {\kern 1pt} {\kern 1pt} {\kern 1pt} {\kern 1pt} {\kern 1pt} {\kern 1pt} {\kern 1pt} {\kern 1pt} {\kern 1pt} {\kern 1pt} {\kern 1pt} {\kern 1pt} {\kern 1pt} {\kern 1pt} {\kern 1pt} {\kern 1pt} if{\kern 1pt} {\kern 1pt} {\kern 1pt} {\kern 1pt} {I_{p\min }} \leqslant x \leqslant {I_{p\max }}\\ 
	&{I_{p\min }}{\kern 1pt} {\kern 1pt} {\kern 1pt} {\kern 1pt} {\kern 1pt} {\kern 1pt} {\kern 1pt} {\kern 1pt} {\kern 1pt} {\kern 1pt} {\kern 1pt} {\kern 1pt} {\kern 1pt} {\kern 1pt} {\kern 1pt} {\kern 1pt} {\kern 1pt} {\kern 1pt} {\kern 1pt} {\kern 1pt} {\kern 1pt} {\kern 1pt} {\kern 1pt} {\kern 1pt} {\kern 1pt} {\kern 1pt} {\kern 1pt} {\kern 1pt} {\kern 1pt} {\kern 1pt} {\kern 1pt} {\kern 1pt} {\kern 1pt} {\kern 1pt} {\kern 1pt} if{\kern 1pt} {\kern 1pt} {\kern 1pt} {\kern 1pt} x \leqslant {I_{p\min }}
	\end{aligned} \right.{\kern 1pt} {\kern 1pt} {\kern 1pt} {\kern 1pt} {\kern 1pt} {\kern 1pt} {\kern 1pt} {\kern 1pt} {\kern 1pt} {\kern 1pt}
	\end{equation}
	The current limit is modeled as follows:
	\begin{enumerate}
		\item Q-priority: $I_{pmax}=\sqrt{I_{max}^2-I_{pcmd}^2}$; if $typeflag=1$ then $I_{pmin}=-I_{pmax}$, else $I_{pmin}=0$.
		\item P-priority: $I_{pmax}=I_{max}$; $I_{pmin}=-I_{max}$;
	\end{enumerate}
	\renewcommand\arraystretch{1.5}
	\begin{table}
		\caption{Parameter definition of DER\_A model \cite{EPRI2019}}
		\label{table1}
		\setlength{\tabcolsep}{8pt}
		\begin{tabular}{|p{30pt}|p{185pt}|}
			\hline
			Parameters& 
			\centerline{Definitions} \\
			\hline
			$T_{rv}$& 
			transducer time constant(s) for voltage measurement\\
			$T_p$& 
			transducer time constant (s) \\
			$T_{iq}$& 
			Q control time constant (s)\\
			$V_{ref0}$& 
			voltage reference set-point $>$ 0 (pu)\\
			$K_{qv}$& 
			proportional voltage control gain (pu/pu)\\
			$T_g$& 
			current control time constant (s)\\
			$Pf_{Flag}$& 
			0 - for constant Q control, and 1 - constant power factor control\\
			$I_{max}$& 
			maximum converter current (pu)\\
			$dbd1$& 
			lower voltage deadband $\leqslant 0$ (pu)\\
			$dbd2$& 
			upper voltage deadband $\geqslant 0$ (pu)\\
			$T_{v}$& 
			time constant on the output of the voltage/frequency cut-off\\
			$V_{lo}$& 
			voltage break-point for low voltage cut-out of inverters\\
			$V_{l1}$& 
			voltage break-point for low voltage cut-out of inverters\\
			$V_{h0}$& 
			voltage break-point for high voltage cut-out of inverters\\
			$V_{h1}$& 
			voltage break-point for high voltage cut-out of inverters\\
			$t_{vl0}$& 
			timer for $V_{l0}$ point\\
			$t_{vl1}$& 
			timer for $V_{l1}$ point\\
			$t_{vh0}$& 
			timer for $V_{h0} $point\\
			$t_{vh1}$& 
			timer for $ V_{h1}$ point\\
			$V_{rfrac}$& 
			fraction of device that recovers after voltage comes back to within $V_{l1} < V < V_{h1}$\\
			$T_{rf}$& 
			transducer time constant(s) for frequency measurement (must be $\geqslant 0.02 s$) \\
			$K_{pg}$& 
			active power control proportional gain\\
			$K_{ig}$& 
			active power control integral gain\\
			$D_{dn}$& 
			frequency control droop gain $\geqslant 0$ (down-side)\\
			$D_{up}$& 
			frequency control droop gain $\geqslant 0$ (up-side)\\
			$f_{emax}$& 
			frequency control maximum error $\geqslant 0$ (pu)\\
			$f_{emin}$& 
			frequency control minimum error $\leqslant 0$ (pu)\\
			$f_{dbd1}$& 
			lower frequency control deadband $\leqslant 0$ (pu)\\
			$f_{dbd2}$& 
			upper frequency control deadband $\geqslant 0$ (pu)\\
			$Freq_{flag}$& 
			0 - frequency control disabled, and 1 - frequency control enabled\\
			$P_{min}$& 
			minimum power (pu)\\
			$P_{max}$& 
			maximum power (pu)\\
			$Tpord$& 
			power order time constant (s)\\
			$dP_{min}$& 
			power ramp rate down $<$ 0 (pu/s)\\
			$dP_{max}$& 
			power ramp rate up $>$ 0 (pu/s)\\
			$V_{tripflag}$& 
			0 $-$ voltage tripping disabled, 1 $-$ voltage tripping enabled\\
			$I_{ql1}$& 
			minimum limit of reactive current injection, p.u.\\
			$I_{qh1}$& 
			maximum limit of reactive current injection, p.u.\\
			$X_{e}$& 
			source impedance reactive $>$ 0 (pu)\\
			$F_{tripflag}$& 
			0 - frequency tripping disabled, 1 - frequency tripping enabled\\
			$PQ_{flag}$& 
			0 - Q priority, 1 - P priority for current limit\\
			$typeflag$& 
			0 - the unit is a generator $I_{pmin} = 0$, 1 - the unit is a storage device and $I_{pmin} = – I_{pmax}$\\
			$V_{pr}$& 
			voltage below which frequency tripping is disabled\\
			\hline
		\end{tabular}
		\label{tab1}
	\end{table}
	\subsection{Static load model}
	In the WECC CMPLDWG, the classic ZIP model is adopted to represent the static load \cite{Wang2019}. The ZIP model consists of constant impedance (Z), constant current (I) and constant power (P) components. It is usually used to represent the relationship between output power and input voltage. The mathematical representation is given as follows.
	\begin{eqnarray}
	\label{Equation.ZIP_P}
	\ {P_{ZIP}} = {P_{0}}\left( {{a_p}{{\left( {\frac{{{V}}}{{{V_0}}}} \right)}^2} + {b_p}\left( {\frac{{{V}}}{{{V_0}}}} \right) + {c_p}} \right)
	\end{eqnarray}
	\begin{eqnarray}
	\label{Equation.ZIP_Q}
	\ {Q_{ZIP}} = {Q_{0}}\left( {{a_q}{{\left( {\frac{{{V}}}{{{V_0}}}} \right)}^2} + {b_q}\left( {\frac{{{V}}}{{{V_0}}}} \right) + {c_q}} \right)
	\end{eqnarray}
	where $P_{ZIP}$ and $Q_{ZIP}$ are active power and reactive power at the bus of interest, $V_0$ is the nominal voltage, $P_{0}$ and $Q_{0}$ are base active/reactive power. $V$ is the voltage magnitude. $a_{p}$, $b_{p}$ and $c_{p}$ are the parameters for active power of the ZIP load, and they satisfy $a_{p}+b_{p}+c_{p}=1$. $a_{q}$, $b_{q}$ and $c_{q}$ are the parameters for reactive power of the ZIP load, and they satisfy $a_{q}+b_{q}+c_{q}=1$. {The first term on the right side of (\ref{Equation.ZIP_P}) represents active power of the constant impedance load, and ${{{P_{0}} \cdot {a_p}} \mathord{\left/ {\vphantom {{{P_{0}} \cdot {a_p}} {V_0^2}}} \right. \kern-\nulldelimiterspace} {V_0^2}}$ is the constant conductance. The second term represents the active power of the constant current load, and ${{{P_{0}} \cdot {b_p}} \mathord{\left/ {\vphantom {{{P_{0}} \cdot {b_p}} {{V_0}}}} \right. \kern-\nulldelimiterspace} {{V_0}}}$ is the constant current. The final term represents the constant power load, and ${P_{0}} \cdot {c_p}$ is the constant power.} 
	\subsection{Electronic load model}
	The electronic load defined in the WECC CMPLDWG is similar to that defined in the software PowerWorld \cite{Wang2019}. The mathematical representation is as follows
	\begin{subequations}
		\label{Equation.Elec}
		\begin{align}
		&  {P_{E,t}} = {c_t} \cdot {P_{E,0}} \quad\quad \label{Elec1}\\[5pt]  
		&  {Q_{E,t}} = {c_t} \cdot {Q_{E,0}}  \quad\quad \label{Elec2}
		\end{align}
	\end{subequations}
	where $P_{E,t}$ and $Q_{E,t}$ are active and reactive power of the electronic load at time $t$, respectively. $c_t$ is a coefficient with respect to the bus voltage, and is defined in Table \ref{table_valuec}. $P_{E,0}$ and $Q_{E,0}$ are base active/reactive power. In Table \ref{table_valuec}, $V_{d1}$ and $V_{d2}$ are two threshold values, and $\alpha$ is a fraction of the electronic load that recovers from low voltage trip. If $\alpha$ is larger than zero, it will be reconnected linearly as the voltage recovers.  ${V_{\min ,t}}$ is a value tracking the lowest voltage but not below $V_{d2}$, and it is a known value at each sample. Its value can be expressed as follows,
	\begin{eqnarray}
	\label{Equation.VminFor}
	\ {V_{\min ,t}} = \max \left\{ {{V_{d2}},\min \left\{ {{V_t},{V_{\min ,t - 1}}} \right\}} \right\}
	\end{eqnarray} 
	The modes depend on the terminal voltage following rules as below:
	\begin{itemize}
		\item If the terminal voltage $V_{t}$ is higher than the threshold value $V_{d1}$, active power and reactive power of the electronic load are constant $P$ and $Q$.   
		\item If the terminal voltage $V_{t}$ is lower than the threshold value $V_{d2}$, active power and reactive power of the electronic load are constant $P$ and $Q$.   
		\item If the voltage $V_{t}$ is between two threshold values $V_{d1}$ and $V_{d2}$ (${V_{d1}} > {V_{d2}}$), active power and reactive power of the electronic load are linearly reduced to zero.
	\end{itemize}
	\begin{table}[]
		\centering
		\caption{Coefficient of Electronic load \cite{Wang2019}}
		\label{table_valuec}
		\renewcommand{\arraystretch}{2.5}
		\setlength{\tabcolsep}{3pt}
		\begin{tabular}{ccc}
			\hline
			Value of $c_t$ & Condition & Mode \\ \hline
			0	&    { { $ \small {{V_t} < {V_{d2}}}$ }}       &   1    \\ \hline
			${\frac{{{V_t} - {V_{d2}}}}{{{V_{d1}} - {V_{d2}}}}}$	&   ${{V_{d2}} \le {V_t} < {V_{d1}},{V_t} \le {V_{\min ,t}}}$         &   2    \\ \hline
			${\frac{{{V_{\min ,t}} - {V_{d2}} + \alpha \cdot ({V_t} - {V_{\min ,t}})}}{{{V_{d1}} - {V_{d2}}}}}$	&     ${{V_{d2}} \le {V_t} < {V_{d1}},{V_t} > {V_{\min ,t}}}$       &   3    \\ \hline
			1	&    ${{V_t} \ge {V_{d1}},{V_{\min ,t}} \ge {V_{d1}}}$        &   4    \\ \hline
			${\frac{{{V_{\min ,t}} - {V_{d2}} + \alpha \cdot ({V_{d1}} - {V_{\min ,t}})}}{{{V_{d1}} - {V_{d2}}}}}$	&   ${{V_t} \ge {V_{d1}},{V_{\min ,t}} < {V_{d1}}}$        &   5   \\ \hline
		\end{tabular}
	\end{table}
	\section{Model validation via simulation}\label{C3}
	In this section, the mathematical model derived in this paper is verified through simulation. The mathematical models of three-phase motor and DER\_A are  tested  on Matlab and PSS/E simultaneously. We compare the performance of the derived mathematical representation with the WECC model embedded in PSS/E to show the accuracy of the derived one. 
	\subsection{Validation of three-phase motors}
\begin{figure}[ht!]
		\centering
		\includegraphics[width=8.5cm]{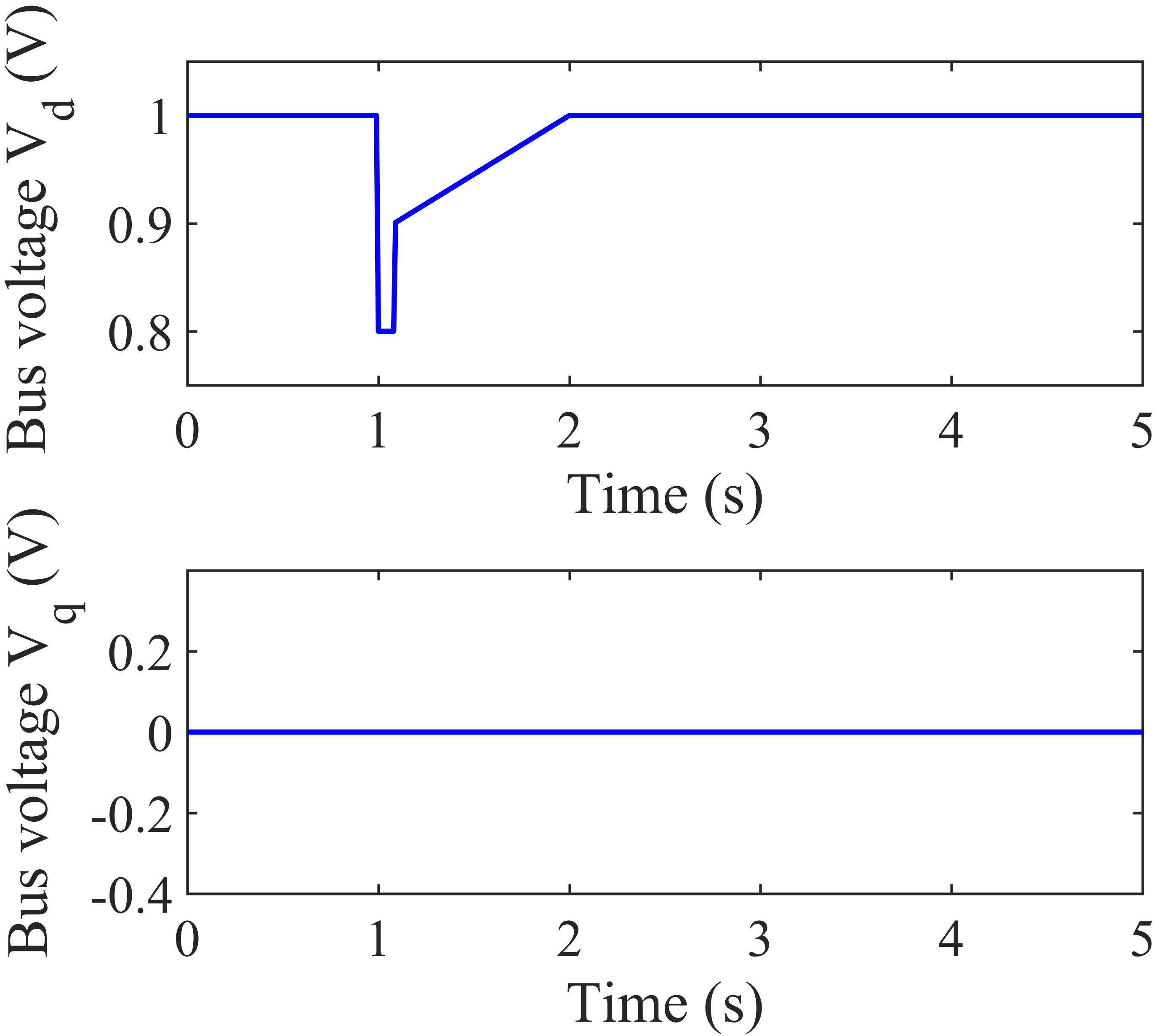}
		\caption{Bus voltages of mathematical and PSS/E model of three-phase motor.}
		\label{V_3phase}
	\end{figure}

\begin{figure}[ht!]
		\centering
		\includegraphics[width=8.5cm]{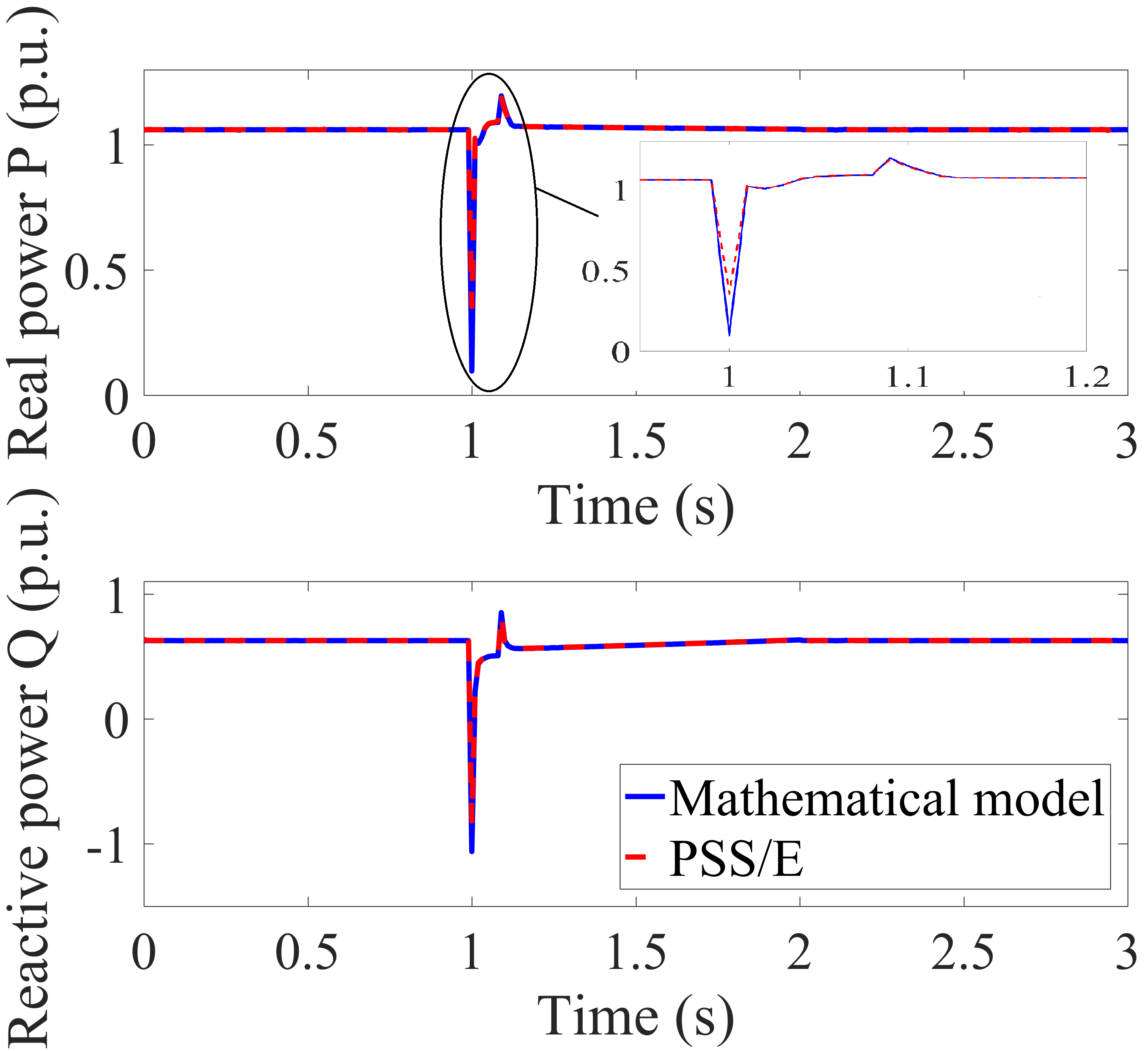}
		\caption{Real and reactive power of mathematical and PSS/E model of three-phase motor A.}
		\label{power_3phaseA}
	\end{figure}

\begin{figure}[ht!]
		\centering
		\includegraphics[width=8.5cm]{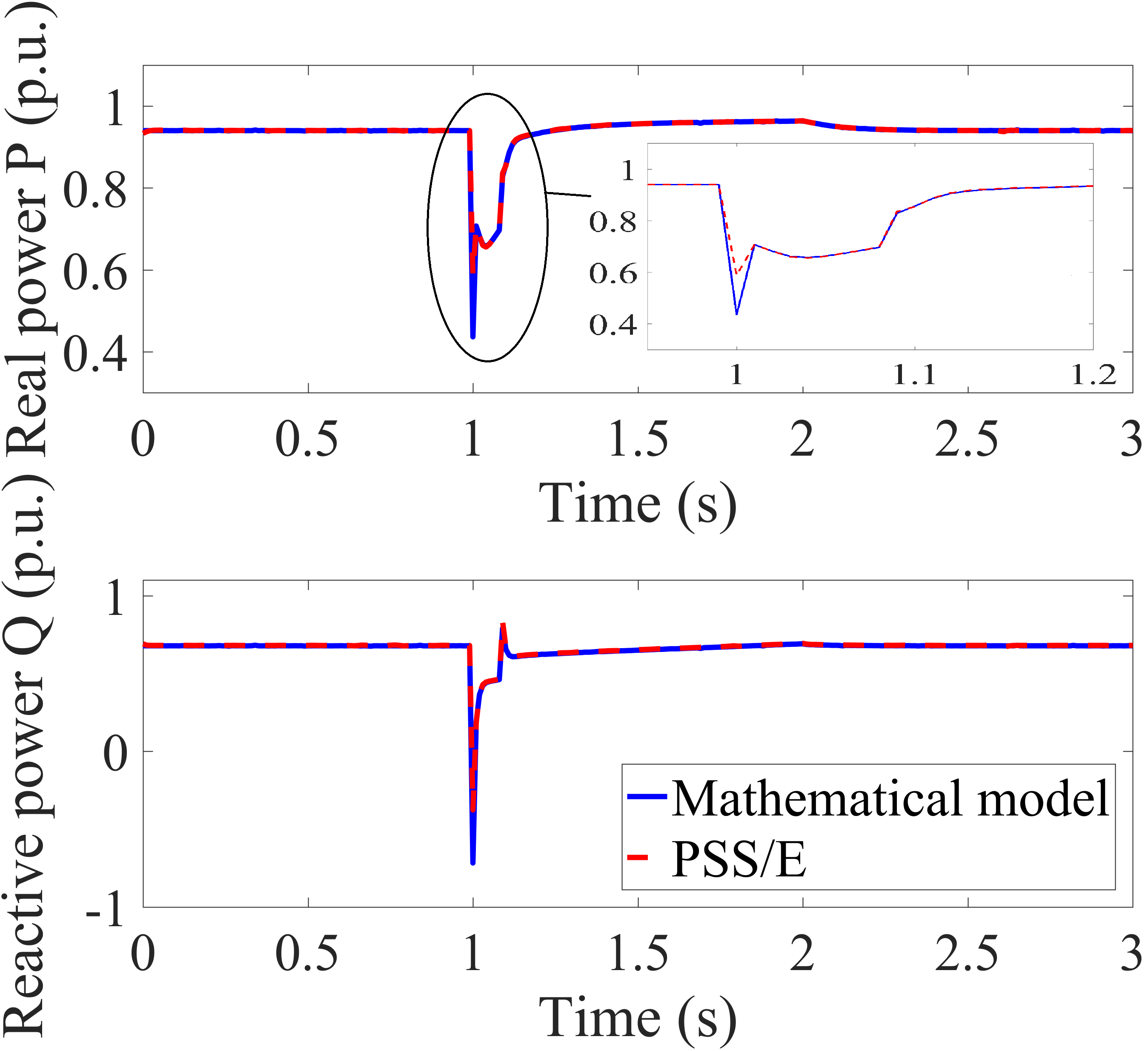}
		\caption{Real and reactive power of mathematical and PSS/E model of three-phase motor B.}
		\label{power_3phaseB}
	\end{figure}

\begin{figure}[ht!]
		\centering
		\includegraphics[width=8.5cm]{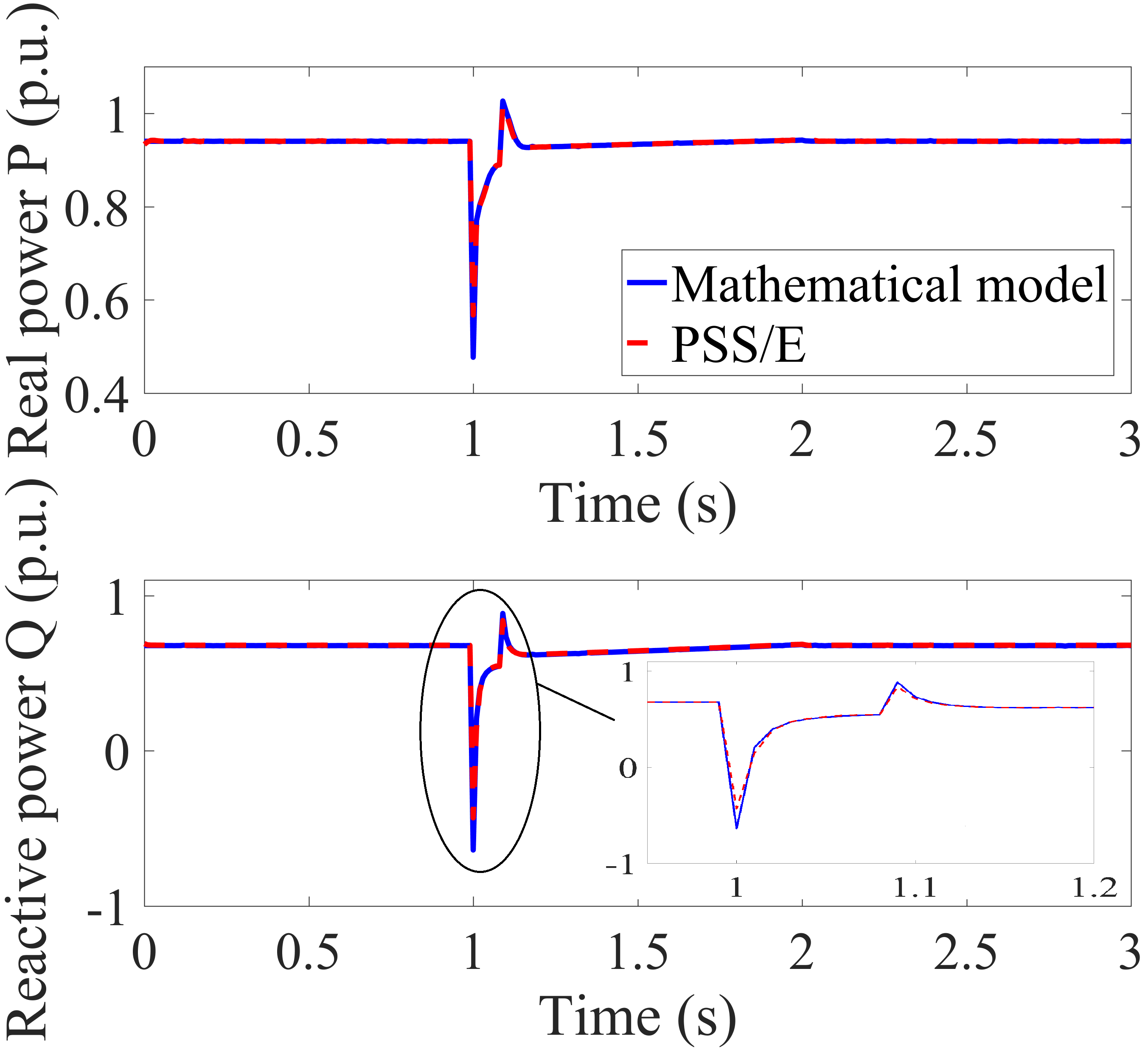}
		\caption{Real and reactive power of mathematical and PSS/E model of three-phase motor C.}
		\label{power_3phaseC}
	\end{figure}

	To verify the proposed model of three-phase motor, we simulated the mathematical model in Matlab and compared it with CMLDBLU2 load model provided by PSS/E. Since here only the mathematical model of three-phase motor is to be validated, the parameters other than three-phase motor in CMLDBLU2 are set to be zero. Moreover, the same bus voltage inputs are given to both models. Consequently, we can compare the output power of the proposed mathematical representation of three-phase motor and that in PSS/E. Refer to \cite{EPRI2019}, the bus voltage input is generated by Equation (\ref{vgen}). The parameters are set as shown in Table \ref{table2}.
	\begin{eqnarray}\label{vgen}
	V\left( t\right) \! =\! \left\{\! \begin{aligned}
	&a{\kern 40pt} if {\kern 5pt} 1\leqslant t\leqslant\left( 1+b/60\right) \!\!\!\!\!\\\!\!\!\!\!\!
	&\!\frac{\left( 1\!-\!d\right) \left(t\!\!-\!\!c\!\!-\!\!1\right)}{b/60-c}\!+\!1 \!{\kern 3pt} if  {\kern 1pt} \left(1\!+\!b/60\right) \leqslant t\leqslant 1\!+\!c\!\!\!\!\!\\\!\!\!\!\!
	&1{\kern 40pt} otherwise
	\end{aligned} \right.
	\end{eqnarray}
	\renewcommand\arraystretch{1.2}
	\begin{table}
		\caption{ParameterS of three-phase motor model \cite{PSERC2017}}
		\label{table2}
		\setlength{\tabcolsep}{9.7pt}
		\begin{tabular}{|c|c|c|c|c|c|}	
			\hline
			\multicolumn{2}{|c|}{Motor A}&\multicolumn{2}{c|}{Motor B}&\multicolumn{2}{c|}{Motor B}\\
			\hline
			$r_{sA}$& 
			0.04&
			$r_{sB}$& 
			0.03&
			$r_{sC}$& 
			0.03\\
			$L_{sA}$& 
			1.8&
			$L_{sB}$& 
			1.8&
			$L_{sC}$& 
			1.8\\
			$L_{pA}$& 
			0.1 s&
			$L_{pB}$& 
			0.16&
			$L_{pC}$& 
			0.16\\
			$L_{ppA}$& 
			0.083&
			$L_{ppB}$& 
			0.12&
			$L_{ppC}$& 
			0.12 \\
			$T_{poA}$& 
			0.092&
			$T_{poB}$& 
			0.1&
			$T_{poC}$& 
			0.1\\
			$T_{ppoA}$& 
			0.002&
			$T_{ppoB}$& 
			0.0026&
			$T_{ppoC}$& 
			0.0026\\
			$H_A$& 
			0.05&
			$H_B$& 
			1&
			$H_C$& 
			0.1\\
			$A_A$& 	0&	$A_B$& 	0&
			$A_C$& 	0\\
			$B_A$& 	0&
			$B_B$& 	0&$B_C$	&0 \\
			$C_A$& 	0&	$C_B$& 	0&
			$C_C$& 	0\\
			$D_A$& 	1&	$D_B$& 	1&
			$D_C$& 	1\\
			$E_{trqA}$& 
			0 &
			$E_{trqB}$& 
			2&
			$E_{trqC}$& 
			2\\
			$p_A$& 	-1&	$p_B$& 	-1&
			$p_C$& 	-1\\
			$q_A$& 	-1&	$q_B$& 	-1&
			$q_C$& 	-1\\
			$\omega_{0A}$& 
			120$\pi$ &
			$\omega_{0B}$& 
			120$\pi$&
			$\omega_{0C}$& 
			120$\pi$\\
			\hline
		\end{tabular}
		\label{tab2}
	\end{table}
	Figure \ref{V_3phase} shows the bus voltage input of three-phase motor. As shown in Figure \ref{power_3phaseA}, \ref{power_3phaseB} and \ref{power_3phaseC} are the dynamic power responses of motor A, motor B and motor C, respectively. The blue solid line denotes the power output of mathematical model, while the red dashed line represents that of CMLDBLU2 in PSS/E. The mean squared errors between the proposed mathematical model and CMLDBLU2 model are shown in Table. \ref{table4}. The small errors show the accuracy of the proposed mathematical model of three-phase model.
	\renewcommand\arraystretch{1.4}
	\begin{table}
		\caption{The mean squared errors between mathematical model and CMLDBLU2 model of three-phase motor.}
		\label{table4}
		\setlength{\tabcolsep}{4.5pt}
		\begin{tabular}{|c|c|c|c|}	
			\hline
			\multirow{2}*{\diagbox{Power}{Motor}}&\multicolumn{3}{c|}{Mean Squared Error (MSE)}\\
			\cline{2-4}
			&Motor A&Motor B&Motor C\\
			\hline
			Real power&$2.1209\times10^{-7}$ &$1.0561\times10^{-5}$  &$1.0762\times10^{-5}$ \\
			\cline{1-4}
			Reactive power&$2.1326\times10^{-5}$ & $1.4474\times10^{-5}$  &		$4.9115\times10^{-5}$ \\
			\hline
		\end{tabular}
		\label{tab4}
	\end{table}
	\subsection{Validation of DER\_A model}
\begin{figure}[ht!]
		\centering
		\includegraphics[width=8.5cm]{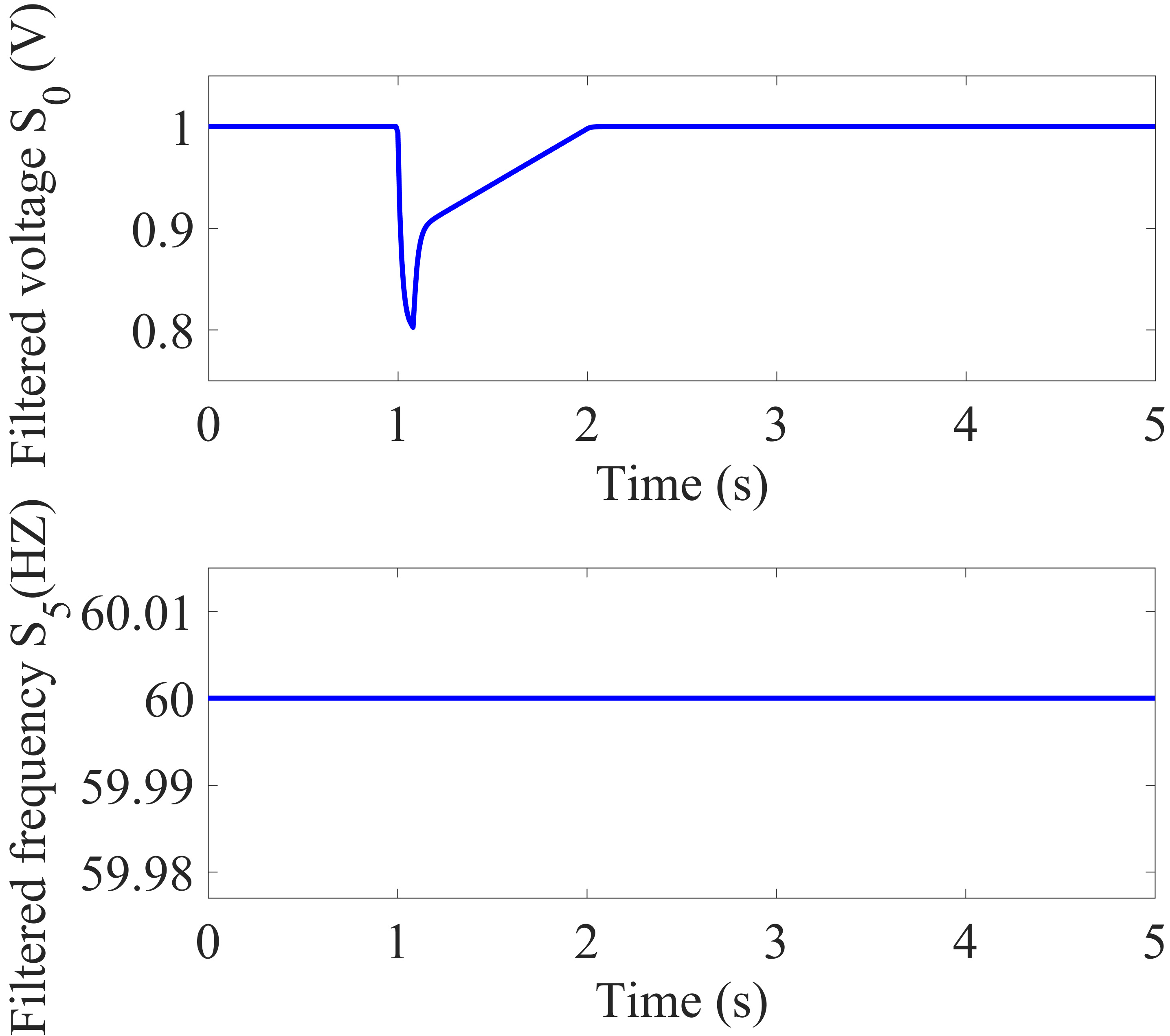}
		\caption{Bus voltages and frequency of mathematical and PSS/E model of DER\_A.}
		\label{V_DERA}
	\end{figure}

\begin{figure}[ht!]
		\centering
		\includegraphics[width=8.5cm]{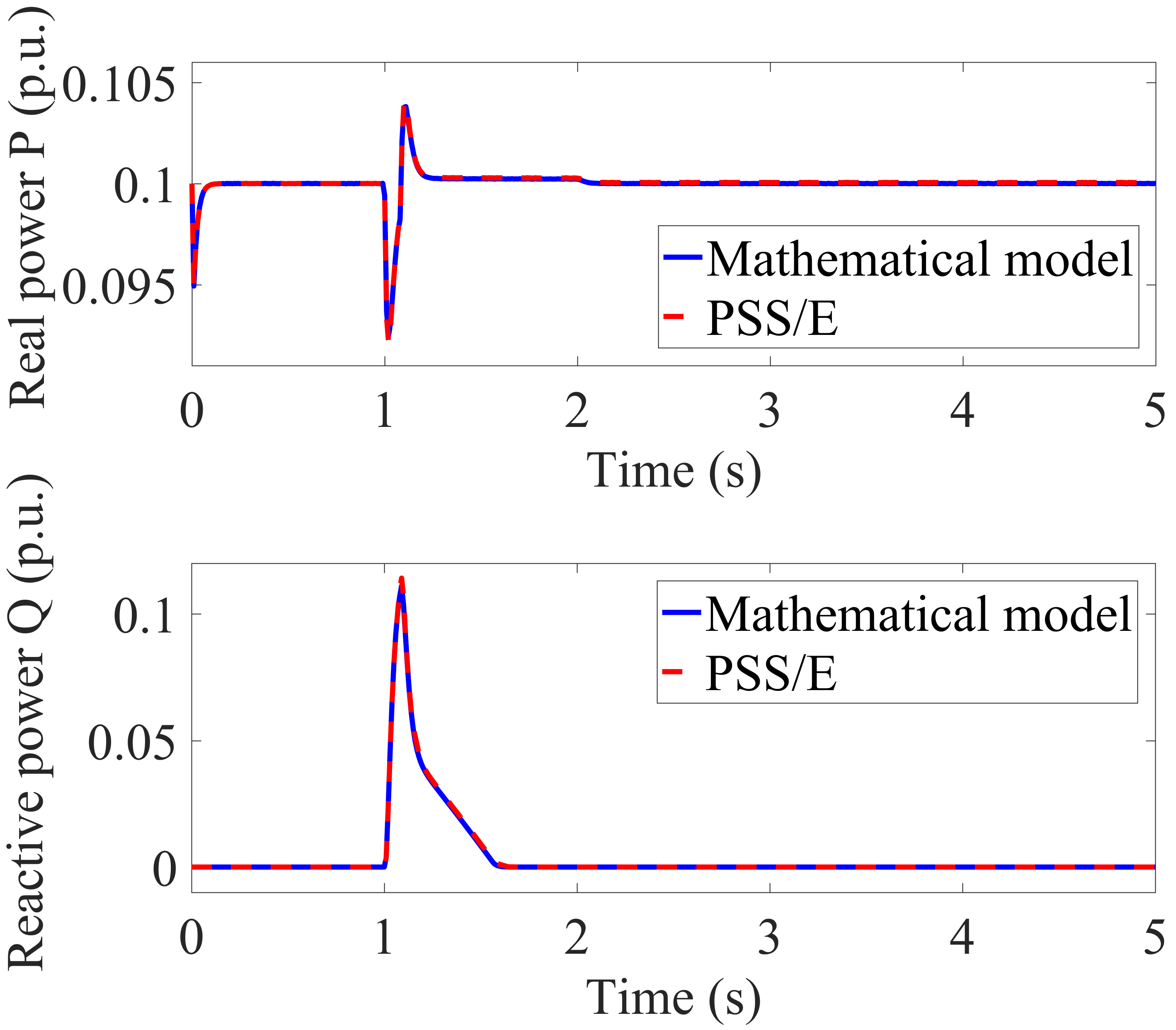}
		\caption{Real and reactive power of mathematical and PSS/E model of DER\_A.}
		\label{power_DERA}
	\end{figure}

	Similar to the verification process of three-phase motor, we simulated the mathematical model of DER\_A in Matlab and adopted DERAU1 provided by PSS/E at the same time. Moreover, the same bus voltage and frequency inputs are given to both models. Consequently, we can compare the output power of the proposed mathematical representation of DER\_A model and that in PSS/E. The voltage input is the same as Equation (\ref{vgen}). The frequency input is set to be 60 HZ. The parameters are set as shown in Table \ref{table3}.
	\renewcommand\arraystretch{1.2}
	\begin{table}
		\caption{Parameter setting of DER\_A model \cite{EPRI2019}}
		\label{table3}
		\setlength{\tabcolsep}{12.7pt}
		\begin{tabular}{|c|c|c|c|}	
			\hline
			Parameters&Values&Parameters&Values\\
			\hline
			$T_{rv}$& 
			0.02 s&
			$T_p$& 
			0.02 s\\
			$T_{iq}$& 
			0.02 s&
			$V_{ref0}$& 
			0 pu\\
			$K_{qv}$& 
			5 pu/pu&
			$T_g$& 
			0.02 s\\
			$Pf_{Flag}$& 
			1&
			$I_{max}$& 
			1.2 pu\\
			$dbd1$& 
			-99 pu&
			$dbd2$& 
			99 pu\\
			$T_{v}$& 
			0.02 s&
			$V_{lo}$& 
			0.44 pu\\
			$V_{l1}$& 
			0.49 pu&
			$V_{h0}$& 
			1.2 pu\\
			$V_{h1}$& 
			1.15 pu&
			$t_{vl0}$& 
			0.16 s\\
			$t_{vl1}$& 
			0.16 s&
			$t_{vh0}$& 
			0.16 s\\
			$t_{vh1}$& 
			0.16 s&
			$V_{rfrac}$& 
			0.7\\
			$T_{rf}$& 
			0.02 s &
			$K_{pg}$& 
			0.1 pu\\
			$K_{ig}$& 
			10 pu&
			$D_{dn}$& 
			20 \\
			$D_{up}$& 
			0 &
			$f_{emax}$& 
			99 pu\\
			$f_{emin}$& 
			-99 pu&
			$f_{dbd1}$& 
			-0.0006\\
			$f_{dbd2}$& 
			0.0006&
			$Freq_{flag}$& 
			0\\
			$P_{min}$& 
			0 pu&
			$P_{max}$& 
			1.1 pu\\
			$Tpord$& 
			0.02 s&
			$dP_{min}$& 
			-0.5 pu/s\\
			$dP_{max}$& 
			0.5 pu/s&
			$V_{tripflag}$& 
			1\\
			$I_{ql1}$& 
			-1 pu&
			$I_{qh1}$& 
			1 pu\\
			$X_{e}$& 
			0.25 pu&
			$F_{tripflag}$& 
			1\\
			$PQ_{flag}$& 
			0&
			$typeflag$& 
			1\\
			$V_{pr}$& 	0.8 pu&	$a$& 	0.8 pu\\
			$b$& 	5&	$c$& 	1 s\\
			\cline{3-4}
			$d$& 	0.9 pu&\multicolumn{2}{c|}{Base: 12.47 kV and 15.0 MVA} \\
			\hline
		\end{tabular}
		\label{tab3}
	\end{table}
	Figure \ref{V_DERA} shows the filtered bus voltage and frequency inputs of DER\_A. Figure \ref{power_DERA} shows is the dynamic power responses of DER\_A. The blue solid line denotes the power output of mathematical model, while the red dashed line represents that of PSS/E. The mean square errors (MSE) of real and reactive power are $1.1053\times10^{-4}$ and $7.3079\times10^{-5}$, respectively. The small error shows the accuracy of the proposed mathematical model of DER\_A.
	\section{Conclusion}\label{C4}
	The WECC composite load model is important for power system monitoring, control and planning, such as stability margin assessment, contingency analysis, assessing the impact of renewable energy, and emergency load control. This paper developed the detailed mathematical model of three-phase motor and DER\_A in WECC composite load model. Several simulations are conducted in matlab and PSS/E. The comparison analysis shows the accuracy of the proposed mathematical representation. This detailed representation is useful for theoretical studies such as stability analysis, parameter identification, order reduction and so forth. 
	
	\bibliographystyle{IEEEtran}
	\bibliography{jiyoulishu}

\end{document}